\documentclass[11pt,a4paper]{article}

\usepackage{psfrag}
\usepackage{amsmath,amssymb,cite}
\usepackage{latexsym}
\usepackage{graphicx}
\usepackage{lscape}
\usepackage{afterpage}
\usepackage[usenames]{color}

\usepackage{abstract}

\usepackage{hyperref}
\hypersetup{
    colorlinks=true,
    linkcolor=blue,
    filecolor=magenta,      
    urlcolor=cyan,
    citecolor=red,
}
\definecolor{kinky}{rgb}{0,.8,.7}
\definecolor{darkgreen}{rgb}{0,.5,.5}
\definecolor{darkorange}{rgb}{.7,.3,0}

\DeclareMathOperator*{\sgn}{sgn}
\DeclareMathOperator*{\diag}{diag}
\DeclareMathOperator*{\E}{\mbox{E}}
\DeclareMathOperator*{\Var}{\mbox{Var}}

\newcommand{\ds}{\displaystyle}
\newcommand{\nexto}{\kern -0.54em}
\newcommand{\dR}{{\rm {I\ \nexto R}}}
\newcommand{\dN}{{\rm {I\ \nexto N}}}
\newcommand{\dZ}{{\cal Z \kern -0.7em Z}}
\newcommand{\dC}{{\rm\hbox{C \kern-0.8em\raise0.2ex\hbox{\vrule
height5.4pt width0.7pt}}}}
\newcommand{\dQ}{{\rm\hbox{Q \kern-0.85em\raise0.25ex\hbox{\vrule
height5.4pt width0.7pt}}}}
\newcommand{\proofbox}{\hspace{\fill}{$\Box$}}

\newtheorem{remark}{Remark}

\newtheorem{algorithm}{Algorithm}

\begin{document}

\author{Authors}

\title{\vspace{-10mm}\bf Sparse Network Optimization for Synchronization}

 \author{
 Regina Sandra Burachik\thanks{Mathematics, UniSA STEM, University of South Australia, Mawson Lakes, S.A. 5095, Australia. E-mail: yalcin.kaya@unisa.edu.au\,.}
 \qquad
Alexander C. Kalloniatis\thanks{Joint and Operations Analysis Division, Defence Science and Technology Group, Canberra 2600, Australia. E-mail: alex.kalloniatis@dst.defence.gov.au\,.}
\qquad
C. Yal{\c c}{\i}n Kaya\thanks{Mathematics, UniSA STEM, University of South Australia, Mawson Lakes, S.A. 5095, Australia. E-mail: yalcin.kaya@unisa.edu.au\,.}
}

\maketitle

\begin{abstract} 
{\vspace{-15mm}\noindent\sf {\bf Abstract.}\ \ We propose new mathematical optimization models for generating sparse dynamical graphs, or networks, that can achieve synchronization. The synchronization phenomenon is studied using the Kuramoto model, defined in terms of the adjacency matrix of the graph and the coupling strength of the network, modelling the so-called coupled oscillators.  Besides sparsity, we aim to obtain graphs which have good connectivity properties, resulting in small coupling strength for synchronization. 
We formulate three mathematical optimization models for this purpose. 
Our first model is a mixed integer optimization problem, subject to ODE constraints, reminiscent of an optimal control problem. As expected, this problem is computationally very challenging, if not impossible, to solve, not only because it involves binary variables but also some of its variables are functions. The second model is a continuous relaxation of the first one, and the third is a discretization of the second, which is computationally tractable by employing standard optimization software. We design dynamical graphs that synchronize, by solving the relaxed problem and applying a practical algorithm for various graph sizes, with randomly generated intrinsic natural frequencies and initial phase variables. We test robustness of these graphs by carrying out numerical simulations with random data and constructing the expected value of the network's order parameter and its variance under this random data, as a guide for assessment.}
\end{abstract}

\begin{verse}
{\em Key words}\/: {\sf Optimization, Sparse graphs, Sparse networks, Synchronization, Kuramoto model, Optimal control, Discretization.}
\end{verse}

\begin{verse} 
{\bf AMS subject classifications.} {\sf Primary 49M20, 49M25\ \ Secondary 92B25}
\end{verse}

\pagestyle{myheadings}
\markboth{}{\sf\scriptsize Sparse Network Optimization for Synchronization\ \ by R. S. Burachik, A. C. Kalloniatis \& C. Y. Kaya}

\section{Introduction} 
\label{intro}

Spontaneous synchronization happens frequently in nature: common examples are fireflies flashing, crickets chirping, planets orbiting and neurons firing. All these phenomena consist of a set of agents that exhibit a cyclic behaviour: these agents are called {\em oscillators} \cite{Forger2017}. Two or more oscillators are said to be {\em coupled} if they influence each other by some physical or chemical process. In our setting, each oscillator is associated with the {\em node} of a {\em graph}, and the whole graph represents a {\em network}. 
The result of such interactions is often {\em synchrony}, which means that all entities end up with the same frequency after a short period of time.  The {\em Kuramoto model} \cite{Kuramoto84} was originally motivated whereby a system of coupled oscillating nodes for sufficient coupling lock on to a common frequency, in this case the mean of the ensemble of natural frequencies of the oscillators, despite the variance in those frequencies.  For a recent literature review of the model see \cite{rodrigues_survey}. 

In the following we summarise salient aspects of the literature relevant to network optimization. The original model, as in \cite{Kuramoto84}, was posed for all-to-all coupling in a complete graph. Since then the pattern of synchronization has been explored for a variety of regular and irregular graphs. 
At one end, sparse regular graphs are the poorest in synchronization, by which we mean they
require high coupling strength to achieve synchronization. For example, synchronization
on regular structures such as trees \cite{DekTay2013} and
rings \cite{RoggAey2004,OchGor2010}
admit a degree of solubility for critical coupling, or bounds thereon.
These are in some respects balanced structures: sparse, with uniform degree distributions.
However, as they scale in size they require increasingly higher values of coupling to achieve synchrony.
Irregular or complex graphs, such as Erd\"os-R\'enyi \cite{Ich2004}, small world \cite{Hong2002}, and the scale-free Barabasi-Alberts \cite{Oh2007} networks synchronize for significantly smaller coupling,
with small world generally the worst of these \cite{Dekk2007}. However, their complexity leads
to very uneven distribution of degrees which may incur some cost, for example imbalances in work or communication loads depending on the physical setting of the synchronization process. In our
case, our interest in such graphs
is for the purpose of synchronized decision-making in organizations 
\cite{Kall2020} where
span of control and cognitive overload
require sparseness and balanced degree
distribution.

In terms of approaches that explicitly use a form of optimization, \cite{Brede2008} employed 
a stochastic hill-climbing approach based on random network link rewirings
and acceptance or rejection on the criterion of improving the Kuramoto order parameter
$r$, to be defined below. Contemporaneously, \cite{TanAoy2008} applied gradient descent
to the order parameter to be able to more efficiently compute optimal weighted networks.
The authors of \cite{KellyGottwald2011} follow an evolution process similar to that of \cite{Brede2008}, but, to avoid hubs in the network, the hill-climbing trajectory only accepts a rewiring that preserves the initial degree distribution. In the same work, the authors propose the use of an additional energy-like function in the selection process.

A Markov Monte Carlo replica exchange approach was developed by \cite{YanMik2012}
to optimise network synchronization against noise, but the method becomes computationally demanding
beyond $N=15$ size graphs. In a second order generalization of the Kuramoto model,
relevant to electrical power systems, \cite{Fazylab2017} used convex optimization to determine optimal placements
of frequencies or link weights within an existing network
structure.
The examples considered in that work were of similarly small size.

A number of optimization approaches have been applied to a related, but
simpler, coupled oscillator model in which the graph Laplacian figures explicitly in the
interaction function, and to which the Kuramoto model approximates under linearization.
Using the master stability function of \cite{PecCar1998}, the authors of \cite{BarPec2002}
observed that graphs with lower Laplacian spectral ratio exhibited improved synchronization.
In this spirit, a stochastic rewiring approach based on improving the spectral ratio was
explored by \cite{Donetti2005,Donetti2006} yielding heterogeneous and typically dense graphs.
Complementary to the spectral ratio approach is that of expander graphs which 
have the property that every partition of the nodes into two
subsets has a number of boundary edges between them that scales with the size of the smallest
subset. Here \cite{Estrada2010} conducted computational search over instances of graphs with such
properties, with results that have excellent expansion and synchronization properties
but that increase in degree as the graphs are scaled in size.

To conclude this summary we mention that \cite{Taylor2020} have provided a construction 
for the Kuramoto model that
yields modified trees -- which are therefore sparse and well-balanced -- that are provably
expanders, with enhanced synchronization. They also showed that a computational search 
over random variations of the links on the leaf nodes of this construction (maintaining
the degree distribution) to find instances with lowest possible spectral ratio 
did not necessarily give graphs of better synchronization.

For more details on similar techniques, which progressively search for a network so as to promote its fitness for synchronization, see \cite[Section 8]{rodrigues_survey}.

In most of these approaches, however, optimization takes the form of a hill-climbing evolution, where, according to certain rules and merit functions,  a sequence of ``mutations" of an initial network, is constructed. After a prescribed number of steps (usually stipulated {\bf a priori}), the last mutation accepted will exhibit some ``optimal" properties, in the sense that the merit functions used in the process have the best values among all the previous mutations. 

In the present paper, we propose, implement and solve a new optimization model for designing graphs that (i) are sparse, (ii) can achieve synchronization in the context of the Kuramoto model \cite{rodrigues_survey} and (iii)  have a relatively narrow range of small degrees. Conditions (i) and (iii) are helpful for ease of implementation of a network, as well as for promoting balanced sharing of
communication- or work-loads.

We mathematically incorporate these three features into our optimization model. Indeed, our optimization variables are the set of all weighted (connected) undirected networks (expressed as $N\times N$ symmetric matrices with nonnegative coefficients), and a set of $N$ functions, namely the phases $\{\theta_i\}$ for $i=1,\ldots,N$ in the Kuramoto system. Our objective function involves a penalty term ensuring sparsity (given by the $\ell_1$ norm of the matrix) and two other terms promoting the small variance of the degrees (see the second and third terms in the objective function in \eqref{model-1}).

To the authors knowledge, there is no available approach that proposes a mathematical optimization model that performs a ``universal search" among all  undirected networks which are able to synchronize by means of the Kuramoto model. 


The initial version of the optimization model we construct 
is an infinite-dimensional mixed-integer programming problem, which is in the form of an optimal control problem.  By relaxing the integer variables, we obtain a continuous-variable model, which is nonsmooth and still infinite-dimensional, and therefore not tractable yet. 
We then discretize the problem, so that optimization software can be employed, as it is done for optimal control problems.  The discretization produces a large-scale problem, which, in general, prohibits any use of nonsmooth solvers.  However, in our particular model, non-smoothness appears in a peculiar form that can be transformed into a smooth, i.e., differentiable, form.  We use tricks and techniques from optimization to convert the problem into a smooth one so that powerful differentiable solvers can be employed to get a solution, even though the problem is still a large-scale and non-convex one.

One should note that a solution to our discrete and smooth optimization model does not represent a graph immediately. A key procedure we
propose, Algorithm~\ref{algo1}, uses solutions of our discrete and smooth optimization model so as to construct increasingly sparse graphs while keeping the spectral ratio of the Laplacian relatively small.  In other words, Algorithm~\ref{algo1} helps us 
\begin{itemize}
    \item[(i)] find a sparse adjacency matrix of the graph such that the system defined by the equations in~\eqref{Kura} below is synchronized, and 
    \item[(ii)] ensure that the resulting graph exhibits a small spectral ratio (defined in~\eqref{Q} below) and 
    \item[(iii)] has a possibly small and narrow range of node degrees. 
\end{itemize}

We present extensive numerical experiments by using this algorithm for graphs with $N=20$ and $N=127$.  We illustrate the fact that, when adding a moderate number of edges, a great improvement in synchronization properties of the graph can be achieved.  We also make comparisons with certain graphs from the literature, in particular from~\cite{Taylor2020}, where graphs are obtained by appending at the bottom level of a hierarchical tree, additional matchings linking the leaf nodes.
The graphs found in this way in \cite{Taylor2020} enhance synchronization compared to
the original tree.
We emphasise the key result: that we can generate graphs with optimal sparseness,
balanced load, and good synchronization properties
starting from a random frequency and
initial phase input.
Because the graphs we generate show low degree heterogeneity, the resultant
synchronization is not so sensitive to
the input frequency choices. 

The paper is organized as follows. In Section \ref{sec1.2}, we recall classical concepts and properties of the Kuramoto model and the spectral ratio of a graph. In Section \ref{Optimization} we introduce our mathematical model that 
promotes synchronization as well as pose its relaxation.  In Section~\ref{smooth}, we discuss discretization of the time-dependent functions in the model and smoothing.  In Section~\ref{num} we present our numerical experiments and provide a detailed interpretation of the numerical results. In Section~\ref{conclusion}, we provide concluding remarks and comment on future work.

\section{Preliminaries}
\label{sec1.2}

In this section we briefly recall the mathematical concepts we will use in our models, including the Tree-expander Construction approach presented in~\cite{Taylor2020}.

\subsection{The Kuramoto Model}
Let $A\in \{0,1\}^{N\times N}$ (i.e., $A$ is an $N\times N$ matrix with all its entries equal to $0$ or $1$) and $\omega = (\omega_1,\ldots,\omega_N)\in \dR^N$. The {\em Kuramoto model} over an $N$-node network can be written as a system of ODEs:
\begin{equation} \label{Kura}
\dot\theta _{i}(t) = \omega_i - \sigma\sum _{j=1}^{{N}}\, A_{ij}\,\sin(\theta _i(t)-\theta _j(t)),\qquad \theta_i(0) = \theta_{i,0}\,,\ \ i=1,\ldots, N,
\end{equation}
where $\theta _i:[0,T]\to \dR$ represents the {\em phase} of the $i$th node, with $\dot\theta _i := d\theta/dt$, and $\theta(t) := (\theta_1(t),\ldots,\theta_N(t))$ is the phase vector in $\dR^N$. Moreover,  $\omega$ is the vector of {\em intrinsic natural frequencies}, $\sigma$ is the {\em coupling strength} of the network, $\theta_0 := (\theta_{i,1},\ldots,\theta_{i,N})$ is the initial phase vector, and $A$ is the {\em adjacency matrix} of a graph with $N$ nodes. 

It is important to distinguish two types of synchronization.
{\it Phase synchronization} refers to asymptotic equality of
the phases, namely the case when $(\theta_i(t)-\theta_j(t))\to 0$ as $t\to\infty,\ \forall \ i,j\in\{1,\ldots,N\}$.  This cannot hold for {\em non-identical oscillators}, where $\omega_i\neq \omega_j$.
{\em Frequency synchronization}, on the other hand, refers to asymptotic equality of the frequencies, namely the case when 
$(\dot{\theta}_i(t) - \dot{\theta}_j(t))\to 0$ as $t\to\infty,\ \forall \ i,j\in\{1,\ldots,N\}$.

Except in certain regular graphs, when coupling is
sufficient for frequency synchronization to set in,
the corresponding phases tend to be clustered
close to each other, namely that
$(\theta_i(t) - \theta_j(t)) \le \varepsilon$ for all $t > t_\ell$ and all $i,j\in\{1,\ldots,N\}$, with some $\varepsilon > 0$ and sufficiently large $t_\ell > 0$, which is referred to as ({\em full}) {\em phase locking}. The higher the coupling is, the smaller the discrepancies between the oscillator phases will be.  In the computational setting that we have, we consider the finite time horizon $[0,T]$ with $T \ge t_\ell$.
This synchronization, however, usually requires dense matrices $A$.  
Viewing links as costly (in resource, in communications,
in load), we seek to keep the topology of the graph as simple as possible. Therefore, it is relevant to find sparse matrices $A$ for which synchronization is achieved.

\subsection{The order parameter}

The so-called {\em order parameter} $r(t)$ is a function of time $t$ defined as
\begin{equation}  \label{r_eqn}
r(t) = \frac{1}{N}\,\left|\sum_{j=1}^N e^{i\,\theta_j(t)} \right|\,,
\end{equation}
where $|\cdot|$ denotes the magnitude of a complex number, with $i = \sqrt{-1}$.  It is not difficult to derive, from~\eqref{r_eqn}, that
\begin{equation}  \label{r2_eqn}
r^2(t) = \frac{1}{N^2}\,\left[N + 2\,\sum_{i<j}^N\,\cos(\theta_i(t)-\theta_j(t)) \right]\,,
\end{equation}
where the usage of $i$ as an index this time should be clear from the context.  As a measure of synchronicity, usually the following long-term average $\overline{r}$ of $r$ is used:
\begin{equation}  \label{rbar_eqn}
\overline{r}(T) = \frac{1}{T-T_{th}}\,\int_{T_{th}}^T r(t)\,dt\,,
\end{equation}
where $T>0$ is a fixed large number and $T_{th}$ is some fraction
of the time-series describing the transient behaviour.

\subsection{Spectral Ratio}
\label{spec_ratio}

Given a graph $G$ with $N$ nodes, call $d_i$ the {\em degree} of node $i$. Collecting all these coordinates in a single vector  $d = (d_1,\ldots,d_N)\in \dN^N$  produces the {\em degree vector} $d$ of $G$. If $A$ is the adjacency matrix of $G$, then the {\em Laplacian} $L$ of $G$ is defined as 
\[
L:=\diag(d)-A\,,
\] 
where $\diag(d)$ is the diagonal matrix with $d_i$ its $i$th diagonal element.  Let $\lambda_2$ and $\lambda_{\max}$ be the second smallest and the largest eigenvalues of $L$, respectively.  Recall that $\lambda_2 \ge 0$. Moreover, if the graph is connected, then $\lambda_2 > 0$. For a connected graph, we denote by 
\begin{equation}\label{Q}
Q_L:=\dfrac{\lambda_{\max}}{\lambda_2}\,,
\end{equation}
the {\em spectral ratio} of the graph $G$.  As alluded
in the introduction, for a somewhat different class of dynamical systems, which the Kuramoto model approximates when linearised, graphs with small $Q_L$ seem to result in highly synchronised graphs, in that the graph is connected and $\sigma$ is small 
\cite{BarPec2002}. For the Kuramoto model such an effect
has been observed in \cite{Taylor2020}.

\subsection{Tree--Expander construction}
\label{tree}

In \cite{Taylor2020}, Taylor, Kalloniatis and Hoek propose a construction procedure for graphs which appends at the bottom level of a hierarchical tree additional matchings linking the leaf nodes. Recall that if the number of levels in a hierarchical tree is $m$, then the total number of nodes in the tree is $p=2^m-1$.  Following the procedure given in \cite{Taylor2020} we have constructed a graph for $m=7$, i.e., $p=127$. We show in Section \ref{127-nodes} a 158-edge graph we have obtained as well as a 158-edge that was obtained in \cite{Taylor2020}, and report their $Q_L$ values. Intuitively, an {\em expander graph} is a graph in which every subset of the vertices that is not ``too large" has a ``large" boundary. Intuitively, this property favours connectivity. The authors show in \cite{Taylor2020} that the graph they construct is an expander graph. 
The resulting graphs have significantly improved
synchronization properties, and are sparse,
of nearly uniform degree distribution and small $Q_L$ (although, for other random leaf node matchings of the tree, the expander construction does not give the smallest possible $Q_L$). 
They therefore
provide a reasonable benchmark to compare
the results of our optimization. In Section~\ref{127-nodes} we compare the properties of these expander graphs with those obtained through our optimization approach described in Algorithm~\ref{algo1}.

\section{Optimization Models}
\label{Optimization}

As mentioned in Section~\ref{intro}, our aims are (i) to find a sparse $A$, and to achieve (ii) synchronization and (iii) a narrow distribution of degrees.  These aims, when realised, are expected to yield a relatively small $Q_L$. To obtain such a matrix $A$, we propose optimization models with decision variables $(A,\theta)$, where $A\in \dR^{N\times N}$ and $\theta$ is the vector function defined in \eqref{Kura}. When the graph is assumed to be undirected, the matrix $A$ is symmetric. Similar models can be posed for the directed case. We fix a random vector $\omega$, and a time $t_c$ after which we impose synchronization among $\theta_i$'s.

\subsection{A Mixed Integer Optimization Problem}

The coefficients of the adjacency matrix are binary, i.e., $A_{ij}\in \{0,1\}$. Let $d$ be the vector of degrees as defined in Section~\ref{spec_ratio}.  Also define
\[
d_{\min}:=\min_{1\le i \le N} d_i\quad\mbox{and}\quad d_{\max}:=\max_{1\le i \le N} d_i.
\]
Fix positive parameters $c_1,c_2,c_3$, $t_c$, as well as a small $\varepsilon > 0$. We propose the following optimization model for the problem that we have.
\begin{equation}  \label{model-1}
\left\{\begin{array}{ll}
& \hbox{\em  objective}: \\
& \hspace*{10mm} \ds\min_{A,\theta}\,\, c_1\sum_{i,j} |A_{i j}| + c_2\,\dfrac{d_{\max}}{d_{\min}} + c_3\,\left(\dfrac{d_{\max}}{d_{\min}}-1\right) \\[6mm]
& \hbox{\em subject to the following constraints.} \\
&\\
& \hbox{\em integrality}:\ \ A_{ij}\in \{0,1\} \ \ \forall i,j,\\
&\\
& \hbox{\em symmetricity}:\ \ A_{ij}=A_{ji} \ \ \forall i,j,\\
&\\
& \hbox{\em oscillator dynamics}: \\
& \ds\dot{\theta}_i(t) = \omega _{i} + \sigma\,\sum_{{j=1}}^{{N}} A_{ij} \sin(\theta_{j}(t) - \theta_{i}(t))\,,\ \forall t\in[0,T]\,,\quad \theta_i(0) = \theta_{i,0}\,,\ \ \forall i\,, \\
&\\
& \hbox{\em synchrony}:\ \ \max_{t\ge t_c}\left\{0, \ds\left(\sum_{ij} |\theta_i(t) -\theta_j(t)|\right)-\varepsilon\right\}=0\,, \ \forall i,j\,.
\end{array}\right.
\end{equation}

The following comments and observations can be made about this model.
\begin{itemize}
\item The first term in the objective function is the $\ell_1$-norm of the elements of the matrix $A$, which, when minimized, is well-known to promote sparsity of $A$.  The two subsequent terms represent the overall discrepancy between the degrees of the nodes.  It turns out that addition of these terms accelerates computation by 10 to 100 times.
\item We have tried specific software for mixed integer optimization, but due to the large number of variables, the above model has so far been intractable.
\item We use a given value of $\sigma$ in the ODE constraints. It is possible to consider $\sigma$ as another decision variable, if smaller values of $\sigma$ are desired.
\end{itemize}

\subsection{A model with relaxation}
\label{model_L1}

Due to the computational complexity, or intractability, of the model in~\eqref{model-1}, we can consider a relaxation of the model in the following way.  First, combine the coefficients of the matrix and the coupling strength $\sigma$ into a single aggregated variable 
\begin{equation}  \label{Aij}
\widetilde{A}_{ij} := \sigma A_{ij}\,,
\end{equation}
and allow $\widetilde{A}_{ij}$ to take any nonnegative real value, i.e. $\widetilde{A}_{ij} \ge 0$.  Next, we also relax the (integer-valued) degree vector by means of the {\em generalized degree vector} $\widetilde{d} = (\widetilde{d}_1,\ldots,\widetilde{d}_N)\in \dR^N$ such that
\begin{equation}\label{deg2}
\widetilde{d}_i:=\sum_{j=1}^N \widetilde{A}_{ij}\,,
\end{equation}
for $i=1,\ldots,N$. Subsequently, the maximum and minimum degree variables we used in the previous model are modified accordingly as
\begin{equation}\label{alfa2}
\widetilde{d}_{\min}:=\min_{1\le i \le N} \widetilde{d}_i\quad\mbox{and}\quad \widetilde{d}_{\max}:=\max_{1\le i \le N} \widetilde{d}_i\,.
\end{equation}
The relaxation \eqref{Aij} allows to consider the ODE \eqref{Kura} with no $\sigma$ appearing on the right-hand side, and, more importantly, to replace the integrality constraints used in the $\ell_1$-model in \eqref{model-1} by non-negativity constraints.  
The use of these aggregated variables allows us to propose a model which can be more easily implemented. So, fix positive parameters $c_1,c_2$ and $c_3$, as before, and consider the following optimization problem.

\begin{equation}\label{model-2}
\left\{\begin{array}{ll}
& \hbox{\em  objective}: \\
& \hspace*{10mm} \ds\min_{\widetilde{A},\theta}\,\, c_1\sum_{i,j} |\widetilde{A}_{i j} |
+ c_2\,\dfrac{\widetilde{d}_{\max}}{\widetilde{d}_{\min}} + c_3\,\left(\dfrac{\widetilde{d}_{\max}}{\widetilde{d}_{\min}}-1\right) \\[6mm]
& \hbox{\em subject to the following constraints.} \\
&\\
& \hbox{\em symmetricity}:\ \ A_{ij}=A_{ji} \ \ \forall i,j, \\
&\\
& \hbox{\em nonnegativity}:\ \ \widetilde{A}_{ij} \ge 0\,, \ \ \forall i,j, \\
&\\
& \hbox{\em oscillator dynamics}: \\
& \ds\dot{\theta}_i(t) = \omega _{i} + \sum_{{j=1}}^{{N}} \widetilde{A}_{ij} \sin(\theta_{j}(t) - \theta_{i}(t))\,,\ \forall t\in[0,T]\,,\quad\theta_i(0) = \theta_{i,0}\,,\ \ \forall i\,, \\
&\\
& \hbox{\em synchrony}:\ \ \max_{t\ge t_c}\left\{0, \left(\sum_{ij} |\theta_i(t) -\theta_j(t)|\right)-\varepsilon\right\}=0\,,\ \ \forall i,j\,.
\end{array}\right. \\
\end{equation}
Since every term in \eqref{deg2} is multiplied by $\sigma$, we have that
\[
\dfrac{\widetilde{d}_{\max}}{\widetilde{d}_{\min}}=\dfrac{{d}_{\max}}{{d}_{\min}}\,.
\]
Therefore, the terms involving $(d_{\max}/d_{\min})$ in the model~\eqref{model-1} remain unchanged by the relaxation, in that $(\widetilde{d}_{\max}/\widetilde{d}_{\min})$ still truly means the ratio of the maximum and minimum degrees in the graph.

\section{Discretization and Smoothing}
\label{smooth}

As it stands in Section~\ref{model_L1}, the $\ell_1$-model is nondifferentiable.   On the other hand, it is well-known that a smooth re-formulation of the $\ell_1$-model can be obtained by using standard nonlinear programming techniques---see e.g.~\cite{NocWri2006,VosMau2006}.   This allows us to use smooth optimization solvers for solving the discretized $\ell_1$-model.

With the elements $\widetilde{A}_{ij}$s of the matrix $\widetilde{A}$ interpreted as constant control functions and $\theta_i$s interpreted as the state variables, the model in \eqref{model-2} is an optimal control problem.  We have discretized the $\ell_1$-model, using the Trapezoidal rule, in a fashion similar to the way optimal control problems are discretized; see, for example, \cite{Kaya2010, KayMar2007}.  Challenging optimal control problems have successfully been solved using direct discretization previously \cite{BanKay2013,KayMau2014,AltKaySch2016,BurKayMaj2014}.

Suppose that the ODE in Model~\eqref{model-2} is defined over the time horizon $[0,T]$, with $T = t_c$.  Let 
\begin{equation}  \label{RHS}
f(\theta_1(t),\ldots,\theta_N(t),\widetilde{A}) := \omega _{i} + \sum _{{j=1}}^{{N}} \widetilde{A}_{ij}\sin(\theta _{j}(t)-\theta _{i}(t)),\quad\forall t\in [0,T]\,.
\end{equation}

Consider a regular partition of $[0,T]$, such that $0 = t_0 < t_1 < \ldots < t_M = T$, with $t_{k+1} := t_k + h$ and $h := T/M$.  Let $\theta_i^k$ be an approximation of $\theta_i(t_k)$, $k = 0,1,\ldots,M-1$.  Then the ODE can be discretized by applying the trapezoidal rule:
\begin{eqnarray}  \label{trapezoidal}
\theta_i^{k+1} &=& \theta_i^k + \frac{h}{2}\left(f(\theta_1^k,\ldots,\theta_N^k,\widetilde{A}) 
+ f(\theta_1^{k+1},\ldots,\theta_N^{k+1},\widetilde{A})\right),
\end{eqnarray}
$k = 0,1,\ldots,M-1$.  After this discretization, and the replacement of the ODE with this discretization, the optimization model in~\eqref{model-2} becomes a finite-dimensional optimization problem.  With the number of optimization variables being $[(N-1)NM/2]$, the discretized problem is large-scale even with a moderate graph size $N$ and a moderate partition size $M$. Moreover, the objective function of the problem is nondifferentiable, prohibiting efficient use of smooth optimization methods and software.

Techniques from the optimization literature can be employed to transform the nonsmooth objective function into a smooth one as follows.  Nonsmoothness of the first term (the $\ell_1$-norm) is caused by the moduli, which can be avoided by defining two new variables, $B_{ij}^1$ and $B_{ij}^2$ such that (see e.g. \cite{VosMau2006})
\begin{equation} \label{smooth_moduli1}
\widetilde{A}_{ij} := B_{ij}^1 - B_{ij}^2\,,\quad B_{ij}^1, B_{ij}^2 \ge 0\,.
\end{equation}
Then one can show that
\begin{equation} \label{smooth_moduli2}
|\widetilde{A}_{ij}| = B_{ij}^1 + B_{ij}^2\,.
\end{equation}
The RHS of the ODEs in \eqref{RHS} then becomes
\begin{equation}  \label{RHS2}
f(\theta_1(t),\ldots,\theta_N(t), B_{ij}^1, B_{ij}^2) := \omega _{i} + \sum _{{j=1}}^{{N}} (B_{ij}^1 - B_{ij}^2)\,\sin(\theta _{j}(t)-\theta _{i}(t))\,,\quad \forall t\in [0,T].
\end{equation}
To tackle the second and third terms in the objective function, one can define the new optimization variables (see e.g. \cite{NocWri2006})
\[
\alpha := \widetilde{d}_{\max}\quad\mbox{and}\quad \beta := \widetilde{d}_{\min}\,.
\]
Then
\[
\min\frac{\widetilde{d}_{\max}}{\widetilde{d}_{\min}} \equiv \min\frac{\alpha}{\beta}
\]
with
\begin{eqnarray}
&& \sum_{j=1}^N \widetilde{A}_{ij} \le \alpha\,,\quad\mbox{for } i = 1,\ldots,N\,, \label{alpha_constr} \\
&& \sum_{j=1}^N \widetilde{A}_{ij} \ge \beta\,,\quad\mbox{for } i = 1,\ldots,N\,. \label{beta_constr} 
\end{eqnarray}
The synchrony constraint in Model~\eqref{model-2} is also nonsmooth, which can be replaced by a smooth version given by
\begin{eqnarray}  \label{synchrony}
\sum_{i\neq j}\left(\theta_i^M -\theta_j^M\right)^2 \le \varepsilon\,.
\end{eqnarray}
Incorporation of \eqref{RHS}--\eqref{synchrony} transforms Model~\eqref{model-2} into the smooth optimization problem below.
\begin{equation}\label{model-3}
\left\{\begin{array}{ll}
& \hbox{\em  objective}: \\
& \hspace*{0mm} \ds\min_{B^1, B^2,\theta^k,\alpha,\beta}\,\, 
c_1\sum_{i,j} (B_{ij}^1 + B_{ij}^2) + c_2\,\dfrac{\alpha}{\beta} + c_3\,\left(\dfrac{\alpha}{\beta} - 1\right) \\[6mm]
& \hbox{\em subject to the following constraints.} \\
&\\
& \hbox{\em symmetricity}:\ \ B_{ij}^1 - B_{ij}^2 = B_{ji}^1 - B_{ji}^2\,,\ \ \forall i,j, \\
& \\
& \hbox{\em nonnegativity}:\ \ B_{ij}^1 - B_{ij}^2 \ge 0\,,\ \ \forall i,j, \\
& \\
&\hbox{\em discretization}: \\[2mm]
& \theta_i^{k+1} = \theta_i^k + \ds\frac{h}{2}\left(f(\theta_1^k,\ldots,\theta_N^k, B_{ij}^1, B_{ij}^2) + f(\theta_1^{k+1},\ldots,\theta_N^{k+1}, B_{ij}^1, B_{ij}^2)\right), \\[2mm]
& \hspace*{70mm} k = 0,1,\ldots,M-1\,,\ \ \theta_i^0 = \theta_{i,0}\,,\ \ \forall i, \\[2mm]
& \hbox{\em synchrony}: \\[1mm]
& \ds\sum_{i\neq j}\left(\theta_i^M -\theta_j^M\right)^2 \le \varepsilon\,, \\
& \\
&\hbox{\em auxiliary}:\ \ B_{ij}^1, B_{ij}^2 \ge 0\,, \forall i,j,\\
& \\
& \hspace*{20mm} \sum_{j=1}^N (B_{ij}^1 - B_{ij}^2) \le \alpha\,,\ \ \forall i\,, \\
& \\
& \hspace*{20mm} \sum_{j=1}^N (B_{ij}^1 - B_{ij}^2) \ge \beta\,,\ \ \forall i\,.
\end{array}\right.
\end{equation}
For the above problem a standard differentiable optimization software can now be used to find a solution.

\section{Numerical Implementation and Experiments}
\label{num}

\subsection{Numerical implementation}
\label{sec:implement}

We consider a numerical implementation of the smooth and finite-dimensional model~\eqref{model-3} derived in Section~\ref{smooth}. We choose the initial values for the vector of intrinsic frequencies $\omega_i$ from the normal distribution with mean $2$ and variance $1/2$, namely $\omega_i\sim {\cal N}(2,1/2)$ and the initial conditions $\theta_i(0)$, $i=1,\ldots,N$, are drawn from the uniform distribution over the interval $[-\pi/2, \pi/2]$, $i = 1,\ldots,N$. 

In solving the discretized and re-formulated optimization model~\eqref{model-3}, we have paired up the optimization modelling language AMPL~\cite{AMPL} and the finite-dimensional optimization software Ipopt~\cite{WacBie2006}.  The AMPL--Ipopt suite was run on a 13-inch 2018 model MacBook Pro, with the operating system macOS Mojave (version 10.14.6), the processor 2.7 GHz Intel Core i7 and the memory 16 GB 2133 MHz LPDDR3.  The Ipopt options {\tt tol=1e-8}, {\tt acceptable\_tol=1e-8} and {\tt linear\_solver=ma57} were provided within an AMPL code.

We have used the following choice of the penalty parameters in model~\eqref{model-3}: $c_1=1$, $c_2=c_3=100$. As mentioned before, the use of the terms with $c_2$ and $c_3$ accelerated convergence by 10--100 times, compared with the case when $c_2 = c_3 = 0$ and the Euclidean distance $\sum_{ij} (\widetilde{d}(i) -\widetilde{d}(j))^2$ is penalized instead.  

A numerical solution of~\eqref{model-3} yields nonnegative real values for $\widetilde{A}_{ij} = B_{ij}^1 - B_{ij}^2$.  It should be noted that $\widetilde{A}_{ij}/\sigma$ for some $\sigma>0$ does not necessarily result in ${A}_{ij}\in\{0,1\}$.  Therefore, we implement the following practical prototype algorithm in constructing a sparse binary adjacency matrix $A$ from the real matrix $\widetilde{A}$.  \\

\noindent
{\bf Prototype Algorithm}\ \\[-7mm]
\begin{description}
\item[Step 1] ({\em Initialization}) Using some $\omega_i\in{\cal N}(2,1/2)$, $\theta_i(0)\in{\cal U}[-\pi/2, \pi/2]$, $i = 1,\ldots,N$, and integer $M>0$, find a solution $\widetilde{A}_{ij}$ to the relaxed model~\eqref{model-3}.  Set $\varepsilon > 0$ small, $\eta_0 = 0$, $T = 2000$ and $t_c = 1500$. Set the coupling strength $\sigma > 0$ large.  Set $k = 0$.
\item[Step 2] ({\em Generate adjacency matrix $A^{(k)}$}) Let\ \ ${A}_{ij}^{(k)} := \sgn(\max\{0,\widetilde{A}_{ij} - \eta_k\})$,\ \ for all $i,j = 1,\ldots,N$.
\item[Step 3] ({\em Solve Kuramoto dynamics}) Solve~\eqref{Kura} over $t\in[0,T]$, with  $\sigma$, $\omega_i$ and $\theta_i(0)$ from Step~1 and $A^{(k)}$ from Step~2.
\item[Step 4] ({\em Test synchrony to potentially increase sparsity of $A^{(k)}$}) If\ \  
$\max_{t_c\le t\le T}\sum_{ij} |\dot\theta_i(t) - \dot\theta_j(t)| \le \varepsilon$\,,\ \ then choose $\eta_{k+1} > \eta_k$, set $k := k+1$ and go to Step 2.
\item[Step 5] ({\em Stopping}) If $k \ge 1$, then declare ``{\em Algorithm successful},'' return the adjacency matrix $A = A^{(k-1)}$.  Otherwise, declare ``{\em Algorithm unsuccessful}.'' Stop.
\end{description}

In Step~2 of the algorithm above, we define the entries ${A}_{ij}^{(k)} := \sgn(\max\{0,\widetilde{A}_{ij} - \eta_k\})$ of $A^{(k)}$ in each iteration $k$.  Note that, with $\eta_k > 0$, some of the entries $\widetilde{A}_{ij} - \eta_k$ may become negative, but after applying the $\max$ operator one gets ${A}_{ij}^{(k)}\in\{0,1\}$.  In Step~4 of the algorithm, first, synchrony is tested with $\eta_0 = 0$, and then $\eta_k$ is increased in each iteration $k$ (i.e., $\eta_{k+1} > \eta_k$) so as to increase the number of $0$-entries in $A^{(k+1)}$, i.e., to obtain a sparser $A^{(k+1)}$. 

Note that the prototype algorithm above is deemed {\em unsuccessful} (in Step~5) only in the case when $k=0$ and the (binary) adjacency matrix $A^{(0)}$ that is generated by Step~2 (with $\eta_0 = 0$) does not produce synchrony.  In the numerical experiments we have performed (those reported in this paper and those not reported), the prototype algorithm (as well as Algorithm~\ref{algo1} given below) was always successful.

Step~3 and the synchrony verification in Step~4 of the prototype algorithm may be replaced by the verification of a connectivity property, as connected graphs do synchronize.  One such property is the strict positivity of the second smallest eigenvalue of the Laplacian of the graph.  Alternatively, one may also replace the same parts of the prototype algorithm by the requirement of a ``suitably small'' value of $Q_L$ in Step~4, so as to have a ``more desirable'' connectivity property of the graph.  Implementation of the latter criterion leads to a slightly different algorithm from the prototype as described below.

\begin{algorithm}  \label{algo1} \rm \ \\[-7mm]
\begin{description}
\item[Step 1] ({\em Initialization}) Using some $\omega_i\in{\cal N}(2,1/2)$, $\theta_i(0)\in{\cal U}[-\pi/2, \pi/2]$, $i = 1,\ldots,N$, and integer $M>0$, find a solution $\widetilde{A}_{ij}$ to the relaxed model~\eqref{model-3}.  Set $\eta_0 > 0$ small. Set $k = 0$.
\item[Step 2] ({\em Generate adjacency matrix $A^{(k)}$}) Same as in Prototype Algorithm.
\item[Step 3] ({\em Potentially increase sparsity of $A^{(k)}$ while keeping spectral ratio small}) While\ \  $Q_L$\ \ is ``small enough,'' choose $\eta_{k+1} > \eta_k$, set $k := k+1$ and go to Step 2.
\item[Step 4] ({\em Stopping}) Same as Step~5 of Prototype Algorithm.
\end{description}
\end{algorithm}

Although $Q_L$ has not been included in the optimization model directly, it is computed in Step~3 of Algorithm~\ref{algo1} (exogenously) for each new graph generated in Step~2.  We have observed in the numerical experiments that, when $k=0$, $Q_L$ is in general small enough, however the matrix is not sufficiently sparse.  As the sparsity of $A$ is increased, $Q_L$ is in general also increased.  Then, if $Q_L$ is not ``small enough'' any longer, we stop and take the latest graph for which $Q_L$ is acceptable.  In other words, we increase $k$ in Algorithm~\ref{algo1} until sparsity and the value of $Q_L$ are both acceptable.

Note that a solution of the relaxed model~\eqref{model-3} depends on the random selection of the vectors $\omega$ and $\theta(0)$.  Moreover, even if we fix the choice for $\omega$ and $\theta(0)$, each run of the AMPL-Ipopt suite with the same $\omega$ and $\theta(0)$ is observed to generate a different (presumably locally optimal) solution.  This is something expected because of the combinatorial nature of the underlying model~\eqref{model-1}.  So, the solution obtained with Algorithm~\ref{algo1}, or Prototype Algorithm, will be different with each $\widetilde{A}$ found by solving model~\eqref{model-3}. Prototype Algorithm is designed to promote sparsity of the graph, and  Algorithm~\ref{algo1} produces graphs with a relatively small spectral ratio. In the numerical experiments we describe in the next section, we have always used Algorithm~\ref{algo1}. This is because we wanted to make sure our resulting graphs had a relatively small spectral ratio, and the matrices produced by it had a high level of sparsity.  

\begin{remark}  \label{rem:vital} \rm
In applications, there might be a need to find out a small set of the so-called {\em vital nodes}~\cite{VitalNodes} that play a critical role, e.g., in propagating information and/or maintaining network connectivity. With this in mind, define a set of links in a connected network to be {\em critical} if their deletion produces a network which is either disconnected, or it has unacceptable connectivity properties such as a large $Q_L$.  In Step~3 of Algorithm~\ref{algo1}, we discard an adjacency matrix if it produces a disconnected graph, or if $Q_L$ is too large. Therefore, Step~3 in Algorithm~\ref{algo1} can ultimately be also used to identify a set of critical links. However, this potential use of Algorithm~\ref{algo1} is not investigated here as it is outside the scope of the current paper.
\end{remark}

\subsection{Numerical experiments}
\label{sec:numexp}

In Subsections~\ref{20-nodes} and \ref{127-nodes}, Problems with 20 and 127 nodes, respectively, are solved. For simplicity in the computational modelling, we set $t_c := T= 20$  in each of these two problems, where $[0,T]$ is the domain, or time horizon, of $\theta$ as in \eqref{Kura}.  This ``terminal" choice of $t_c$ indeed produces the desired synchronization, way earlier. We also use a moderate grid size, $M=10$, for the time horizon so as to make computations faster.

The reason we have chosen $N=127$ is that this number is of the form $N = 2^m -1$, with $m = 7$, and therefore the 127-node graphs can be compared with those obtained in \cite{Taylor2020}.  In the latter reference, graphs with good connectivity properties are constructed by combining a binary tree with an expander graph, as also briefly explained in Section~\ref{tree}. Our aim in choosing $N=127$ is simply to be able to compare the connectivity properties of constructions such as those in \cite{Taylor2020} with those we obtain through Algorithm~\ref{algo1}.

For each of the 20- and 127-node problems, we proceed as follows. First we obtain an acceptable (binary) adjacency matrix $A$ as an output of Algorithm~\ref{algo1}.  We randomly generate 30 sets of the vectors $\omega$ and $\theta(0)$, whose components are taken from the distributions ${\cal N}(2,1/2)$ and ${\cal U}[-\pi/2, \pi/2]$, respectively.  For each of the 30 randomly generated sets of $\omega$ and $\theta(0)$, we solve the Kuramoto system for a sequence of (fixed) values of $\sigma$.  For this purpose we have used {\sc Matlab}, Release 2019b.  In order to solve~\eqref{Kura}, we implemented the solver {\tt ode15s}, effective for stiff ODEs, with the absolute and relative tolerances of $10^{-6}$.  

Now, with the discrete solution of the phase vector $\theta(t)$ for a fixed value of $\sigma$ at hand, Equations~\eqref{r2_eqn}--\eqref{rbar_eqn} are used to obtain the corresponding $\bar r$.  For numerical integration we have used Simpson's rule.  Then using the 30 different sequences of $\bar r$  obtained for each of the 30 random data, the expected value $\E[\bar r]$ and the variance $\Var[\bar r]$ of $\bar r$ are computed and plotted against the coupling strength $\sigma$.  We also provide for each graph a histogram of the degrees of the nodes.  

The variance plot can be used as a measure of robustness of the connectivity of the graph for $\omega$ and $\theta(0)$ coming from distributions that were also chosen in Algorithm~\ref{algo1}.  The plot provides a quick picture of what coupling strength values a network will be resilient, i.e. will be insensitive, to possible changes in $\omega$ and $\theta(0)$.

The adjacency matrices of the graphs mentioned in what follows can be accessed at the URL in~\cite{adjacency_matrices}.

\begin{figure}[t!]
\begin{center}
\hspace*{4mm}
\includegraphics[width=120mm]{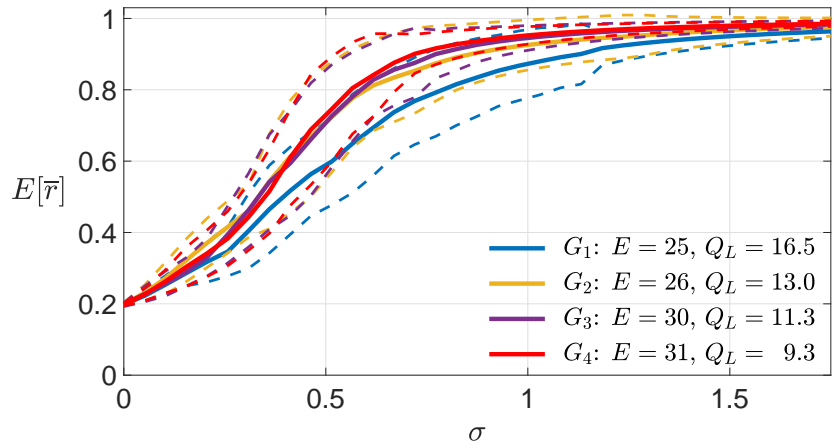} \\
(a) $\E[\bar r]$ vs $\sigma$.
\end{center}
\begin{center}
\hspace*{4mm}
\includegraphics[width=120mm]{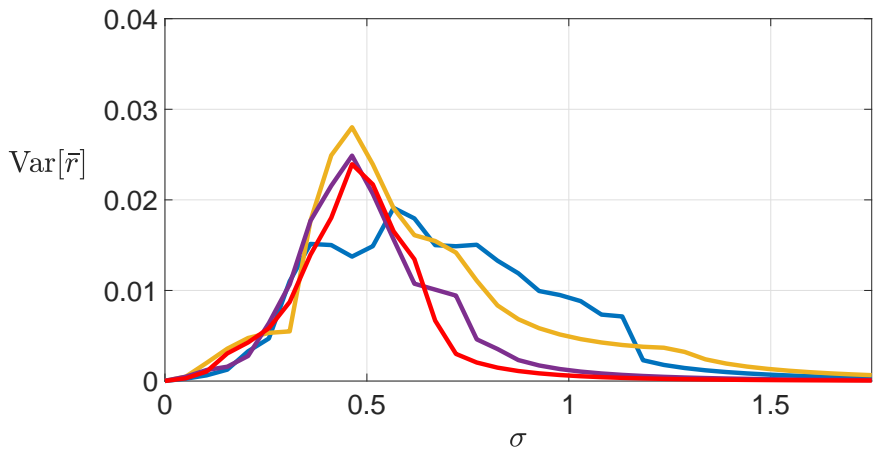} \\
(b) $\Var[\bar r]$ vs $\sigma$.
\end{center}
\caption{\sf $N=20$: Mean and variance of the averaged order parameters of the graphs $G_1$--$G_4$ obtained by Algorithm~\ref{algo1}.} 
\label{fig:V20_combined}
\end{figure}

\subsubsection{20-node networks}
\label{20-nodes}

We used Algorithm~\ref{algo1} for $N=20$. We have chosen some representative instances that illustrate the effectiveness of our approach. Figure~\ref{fig:V20_combined} corresponds to the graphs, denoted $G_i$, $i = 1,\ldots,4$, that we obtained by using Algorithm~\ref{algo1}.

\begin{figure}[t!]
\begin{center}
$G_1$:
\begin{minipage}{65mm}
\includegraphics[width=60mm]{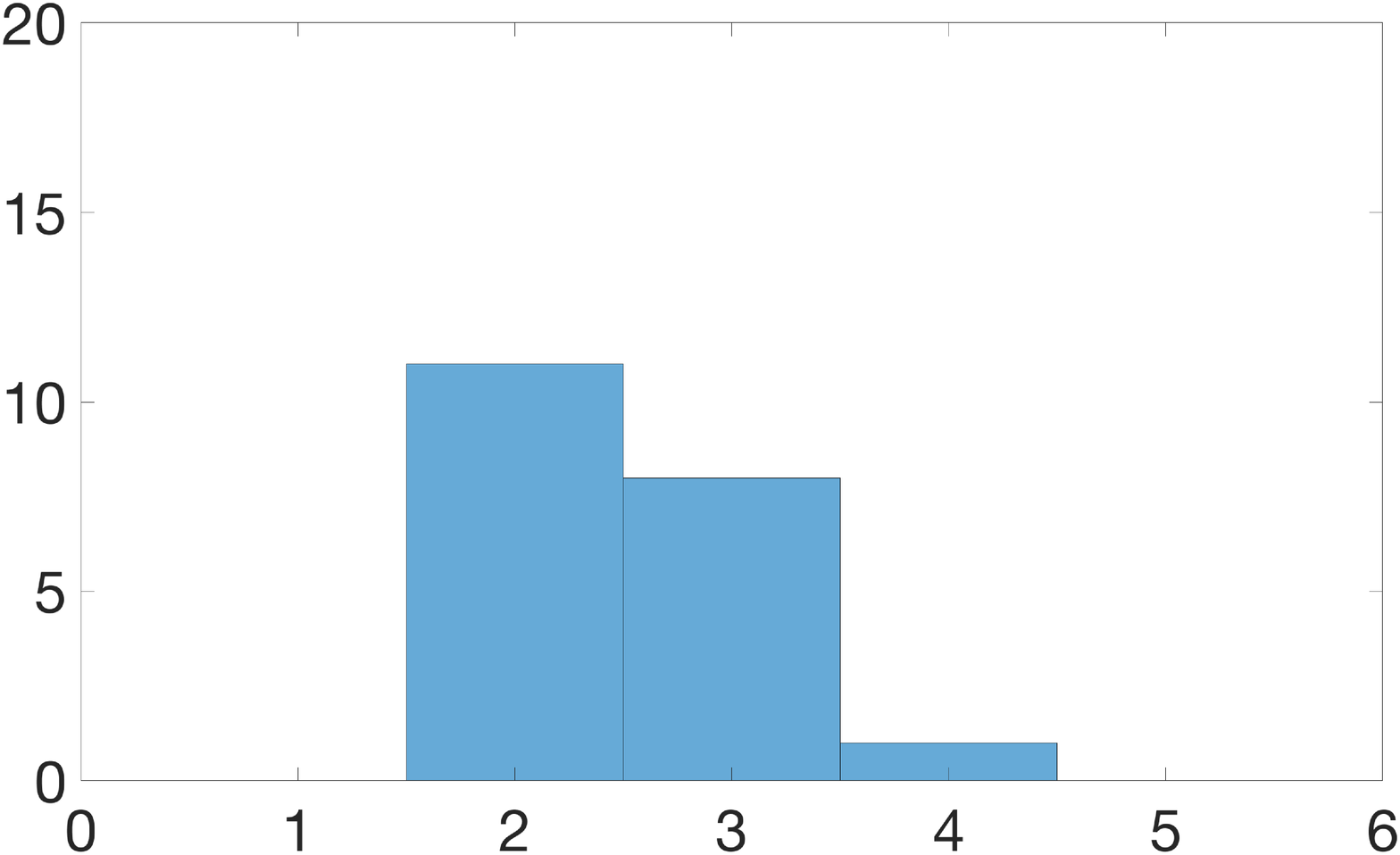}
\end{minipage}
$G_2$:
\begin{minipage}{65mm}
\includegraphics[width=60mm]{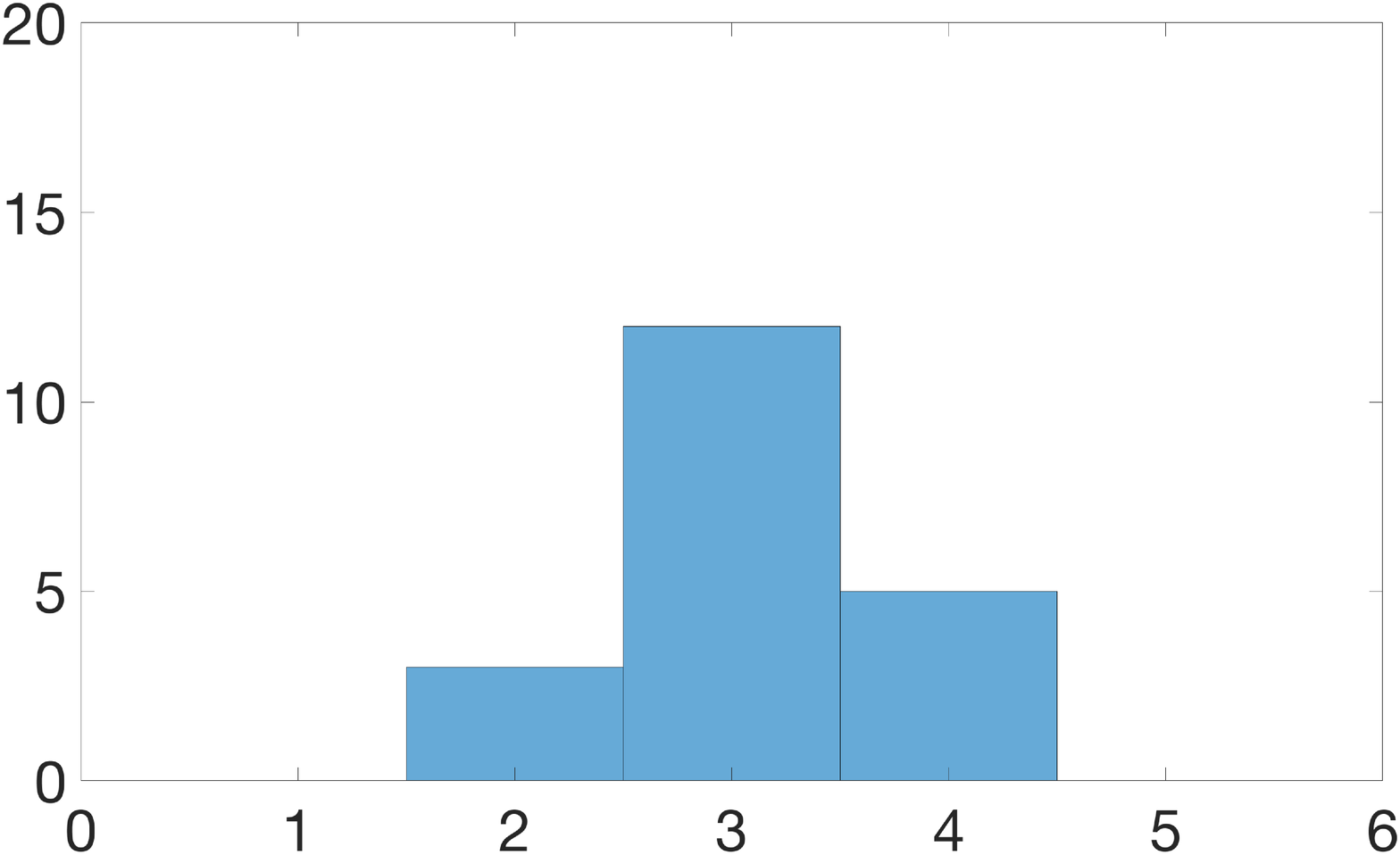}
\end{minipage} \\
$G_3$:
\begin{minipage}{65mm}
\includegraphics[width=60mm]{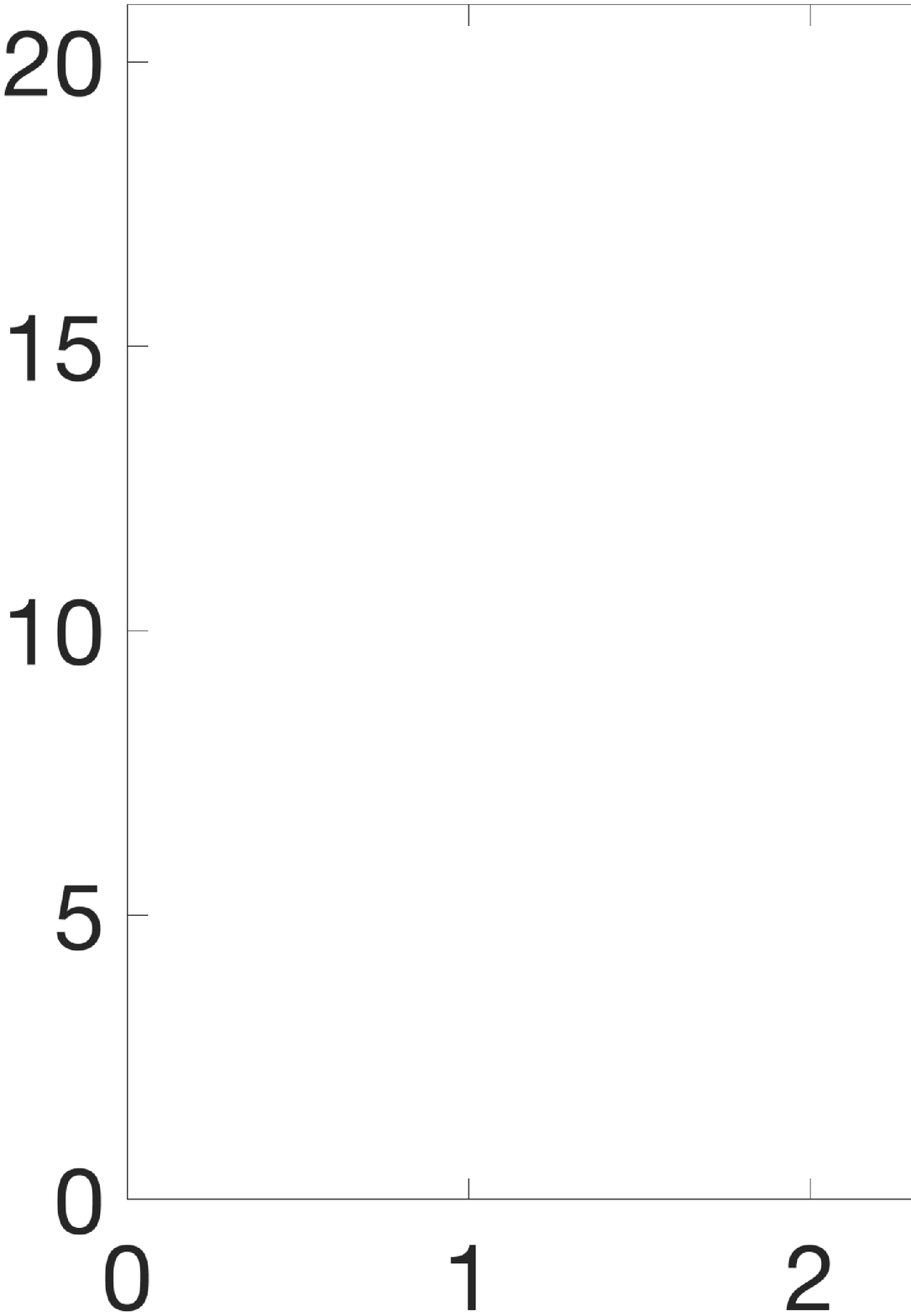}
\end{minipage}
$G_4$:
\begin{minipage}{65mm}
\includegraphics[width=60mm]{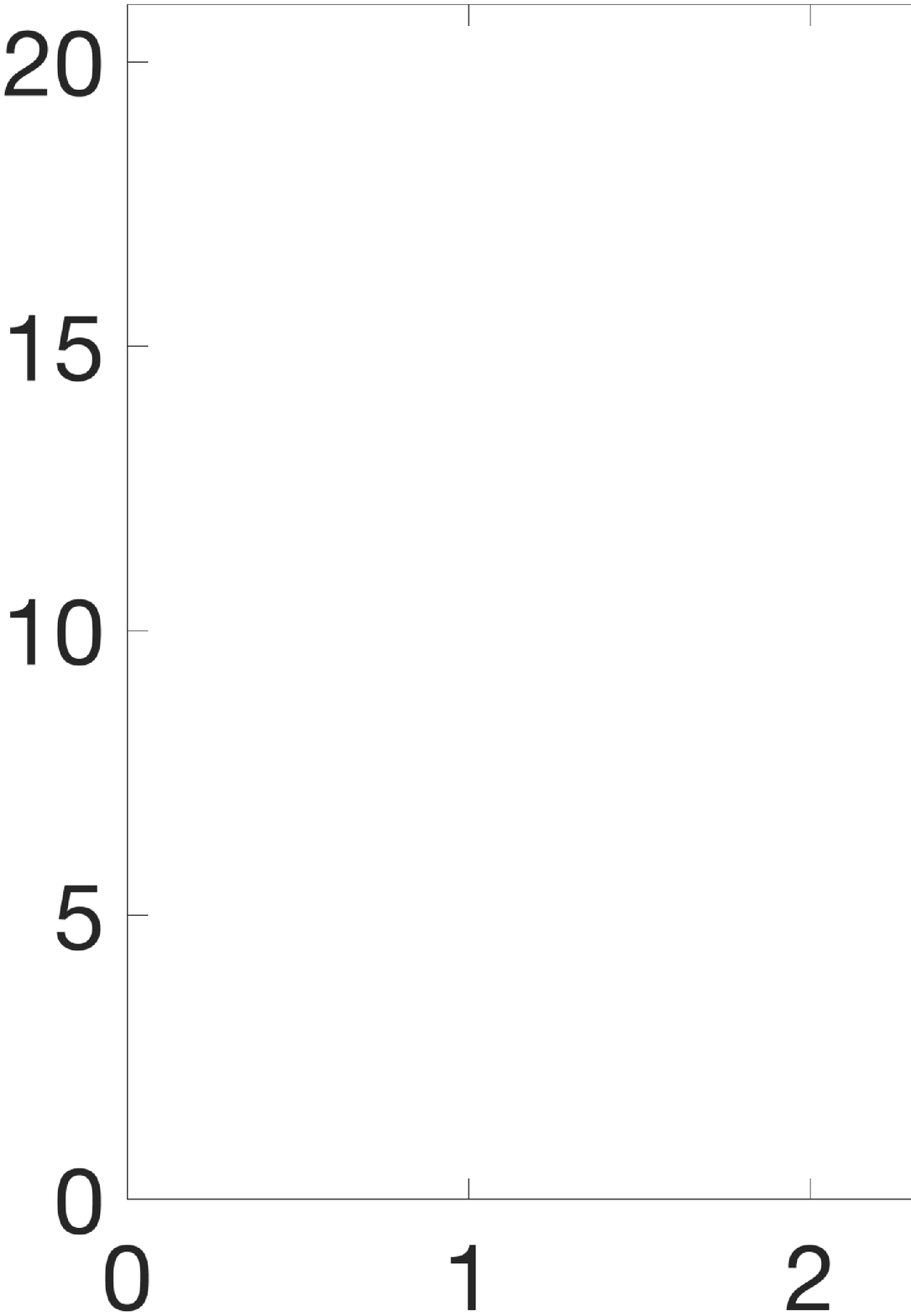}
\end{minipage}
\end{center}
\caption{\sf $N=20$: Degree distributions of the graphs $G_1$--$G_4$ obtained by Algorithm~\ref{algo1}.} 
\label{fig:V20_A_deg_combined}
\end{figure}

In Figure~\ref{fig:V20_combined}(a) we display the plots of the expected values $\E[\bar r]$ of $\bar r$ vs the coupling strength $\sigma$, shown by solid curves which are colour-coded for each $G_i$, $i = 1,\ldots,4$. The dashed curves, also colour-coded, represent the standard deviation from each of the $\E[\bar r]$ curves.  Figure~\ref{fig:V20_combined}(b), on the other hand, depicts the variance $\Var[\bar r]$ of $\bar r$, colour-coded in the same way as in Figure~\ref{fig:V20_combined}(a).  These plots indicate for which values of the coupling strength one can achieve connectivity, and thus synchrony, with a degree of robustness.  Figure~\ref{fig:V20_A_deg_combined} displays the degree histogram of each graph.  

Figure~\ref{fig:V20_combined}(a) and \ref{fig:V20_combined}(b) suggest that the graph $G_1$, with 25 edges and $Q_L = 16.5$, would likely achieve a robust synchrony for $\sigma \ge 1.2$, since, for these values, $\E[\bar r]$ is near $1.0$ and $\Var[\bar r]$ is relatively small.  The graph $G_2$ has just one more edge than $G_1$: With 26 edges, it has a better value of $Q_L$, namely $Q_L=13$.  Figure~\ref{fig:V20_A_deg_combined}(a) shows that its $\E[\bar r]$ curve comes near the value $1.0$ more rapidly than that of $G_1$.  Moreover, one can observe in Figure~\ref{fig:V20_A_deg_combined}(b) that its $\Var[\bar r]$ curve decays relatively more quickly.  As can be seen in Figure~\ref{fig:V20_A_deg_combined}, the degrees of the nodes of both graphs $G_1$ and $G_2$ range between 2 and 4, but a great majority of the nodes of $G_2$ are of degree 3 or 4.

Figure~\ref{fig:V20_A_deg_combined} also depicts that $G_3$ is a {\em cubic} graph, i.e. each of its nodes has degree 3.  Cubic graphs are desirable in many practical applications because of their balanced
workload sharing nature.  We note that $G_3$ has 30 edges but $Q_L=11.3$ is significantly smaller than those for $G1$ and $G2$.  In the plots in Figure~\ref{fig:V20_combined}, we observe that $\E[\bar r]$ approaches $1.0$ and $\Var[\bar r]$ (after making the usual peak) dies down much more rapidly than the those curves for $G1$ and $G2$, pointing to far more robustness in its synchronization.

The graph $G_4$ has 31 edges, and an even better value of $Q_L=9.3$.  The node degrees in this case range from 3 to 4, which can also be regarded as sharing
workloads more evenly compared to $G_1$ and $G_2$.  The ``performance indicators'' of the graph given in Figures~\ref{fig:V20_combined} are even more impressive than before, in terms of the criteria discussed for the previous graphs.  Robustness of its synchronization is more evident with its $\Var[\bar r]$ being much smaller than those for $G_1$, $G_2$ and $G_3$ with larger (synchronizing) values of $\sigma$.

\begin{figure}[t!]
\begin{center}
\hspace*{4mm}
\includegraphics[width=120mm]{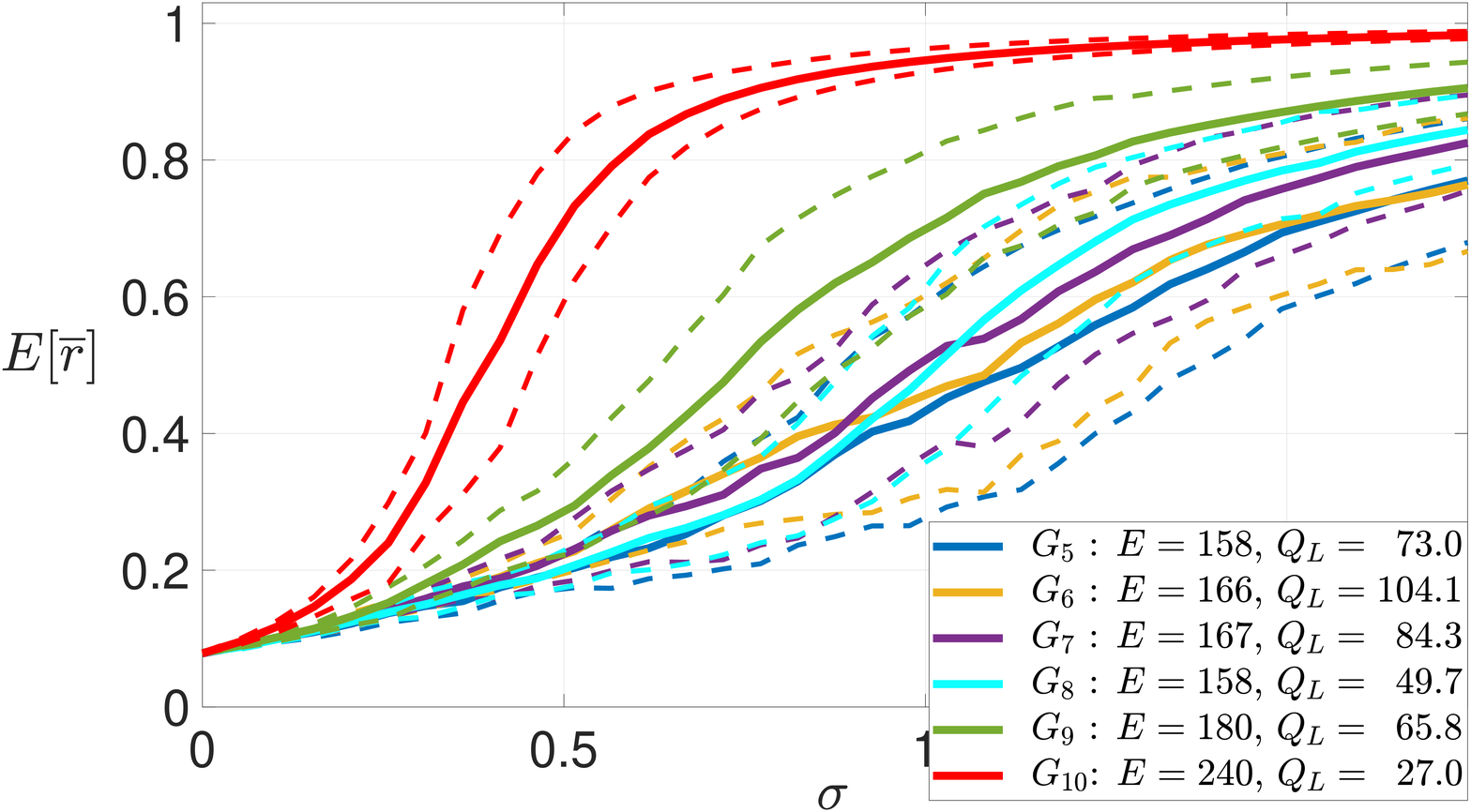} \\
(a) $\E[\bar r]$ vs $\sigma$.
\end{center}
\begin{center}
\hspace*{4mm}
\includegraphics[width=120mm]{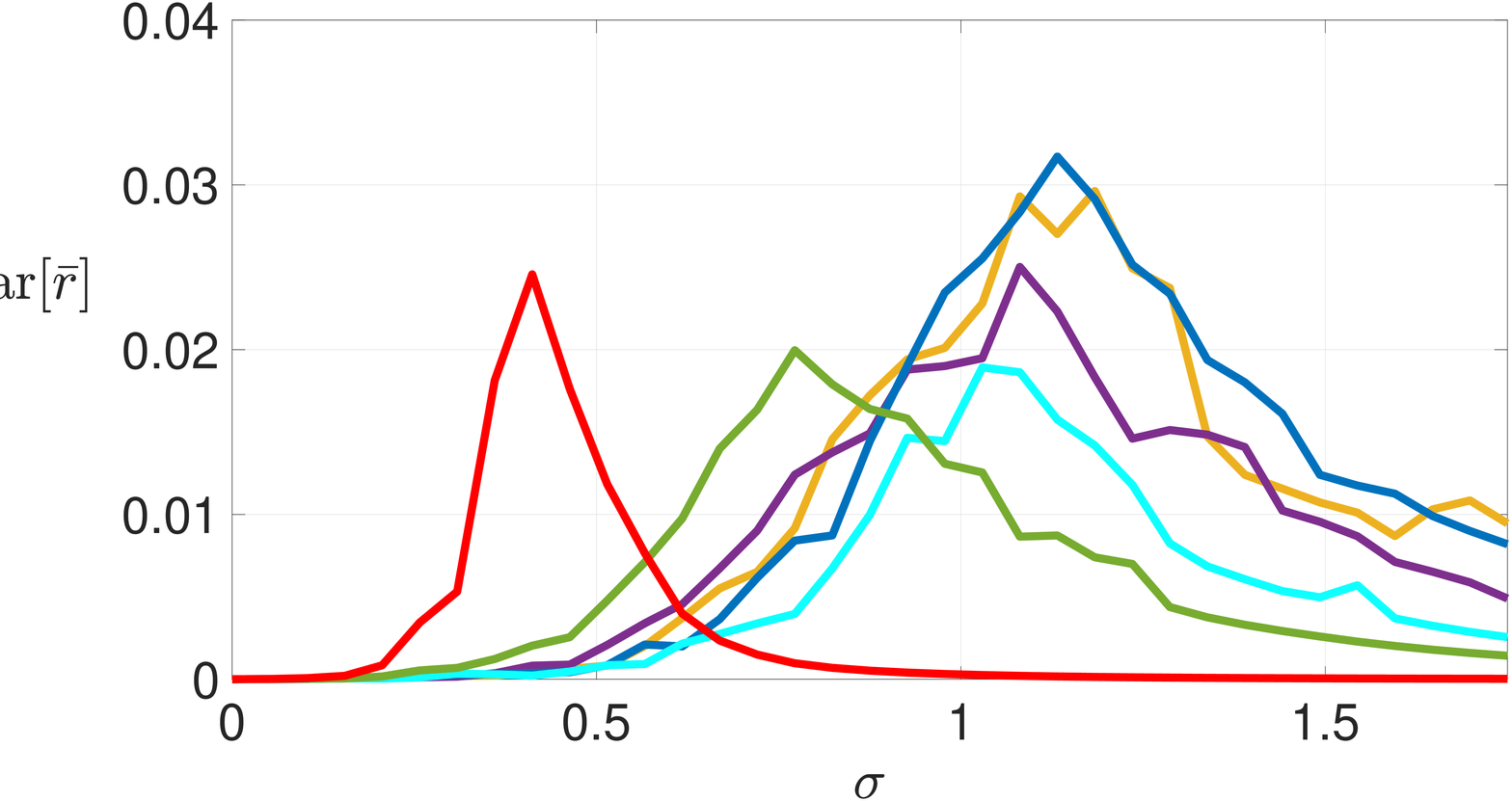} \\
(b) $\Var[\bar r]$ vs $\sigma$.
\end{center}
\caption{\sf $N=127$: Mean and variance of the averaged order parameters of the graphs $G_5$--$G_{10}$ obtained by Algorithm~\ref{algo1}.} 
\label{fig:V127_combined}
\end{figure}

\begin{figure}[t!]
$G_5$:
\begin{minipage}{65mm}
\includegraphics[width=60mm]{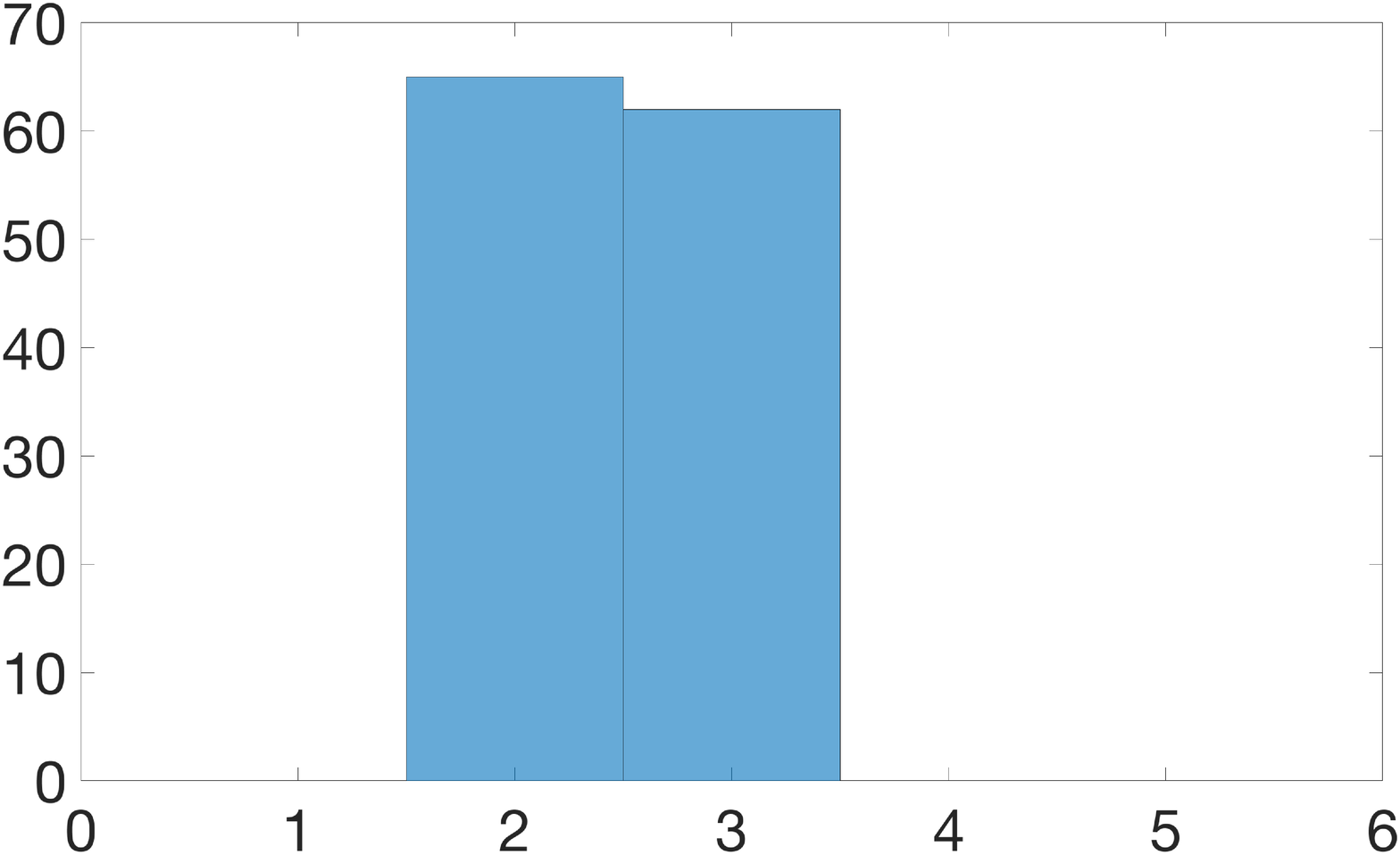}
\end{minipage}
$G_6$:
\begin{minipage}{65mm}
\includegraphics[width=60mm]{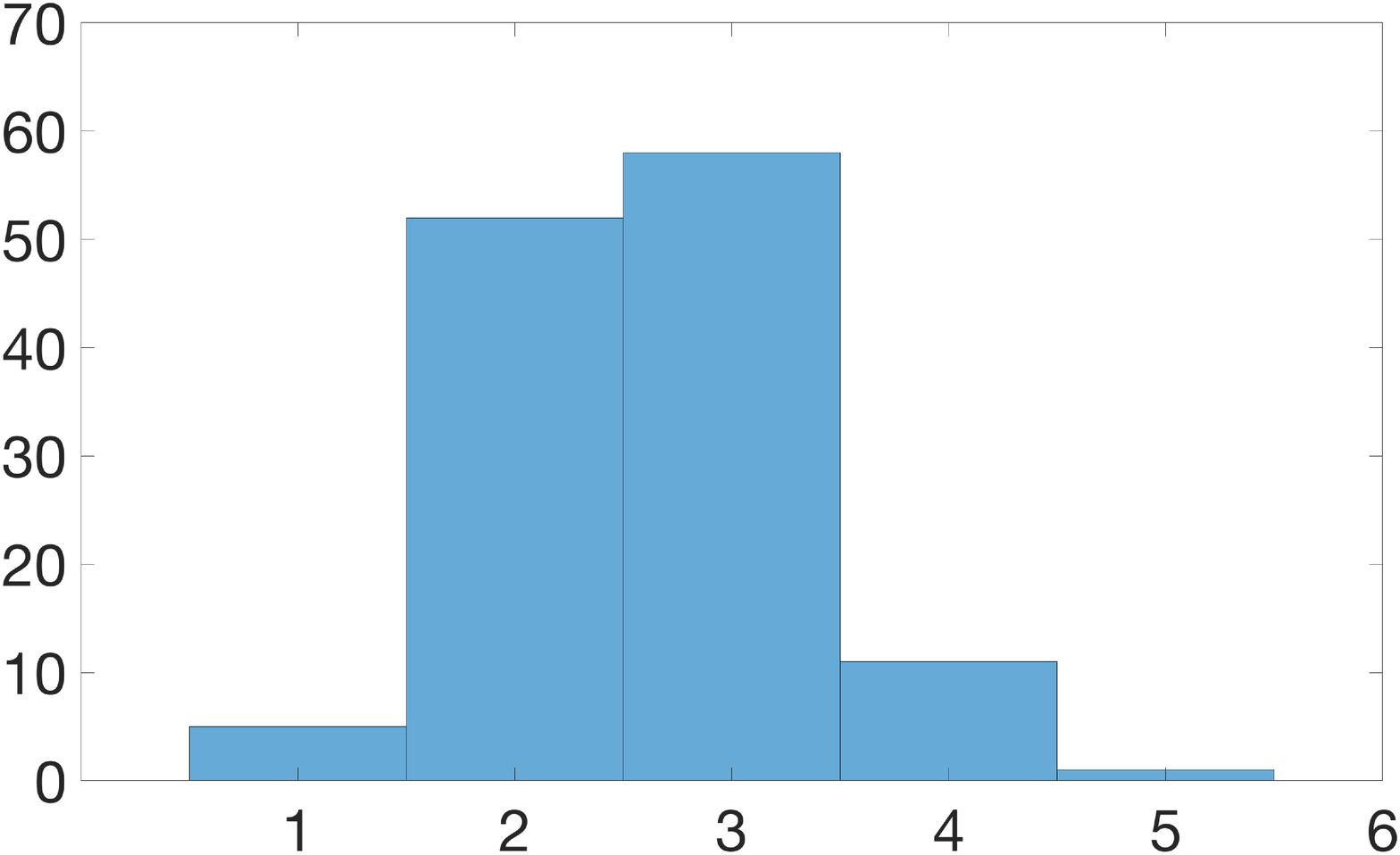}
\end{minipage} \\
$G_7$:
\begin{minipage}{65mm}
\includegraphics[width=60mm]{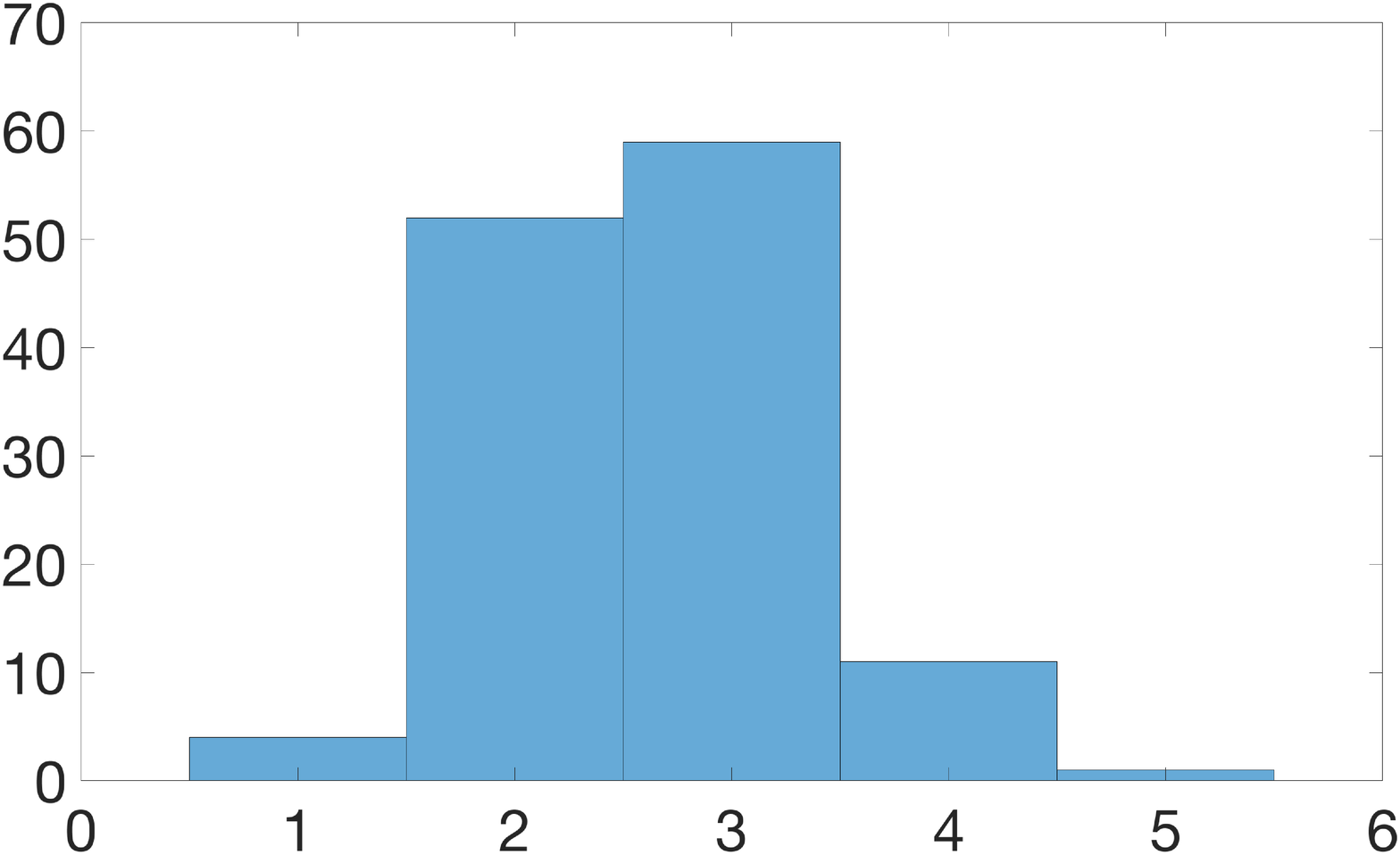}
\end{minipage}
$G_8$:
\begin{minipage}{65mm}
\includegraphics[width=60mm]{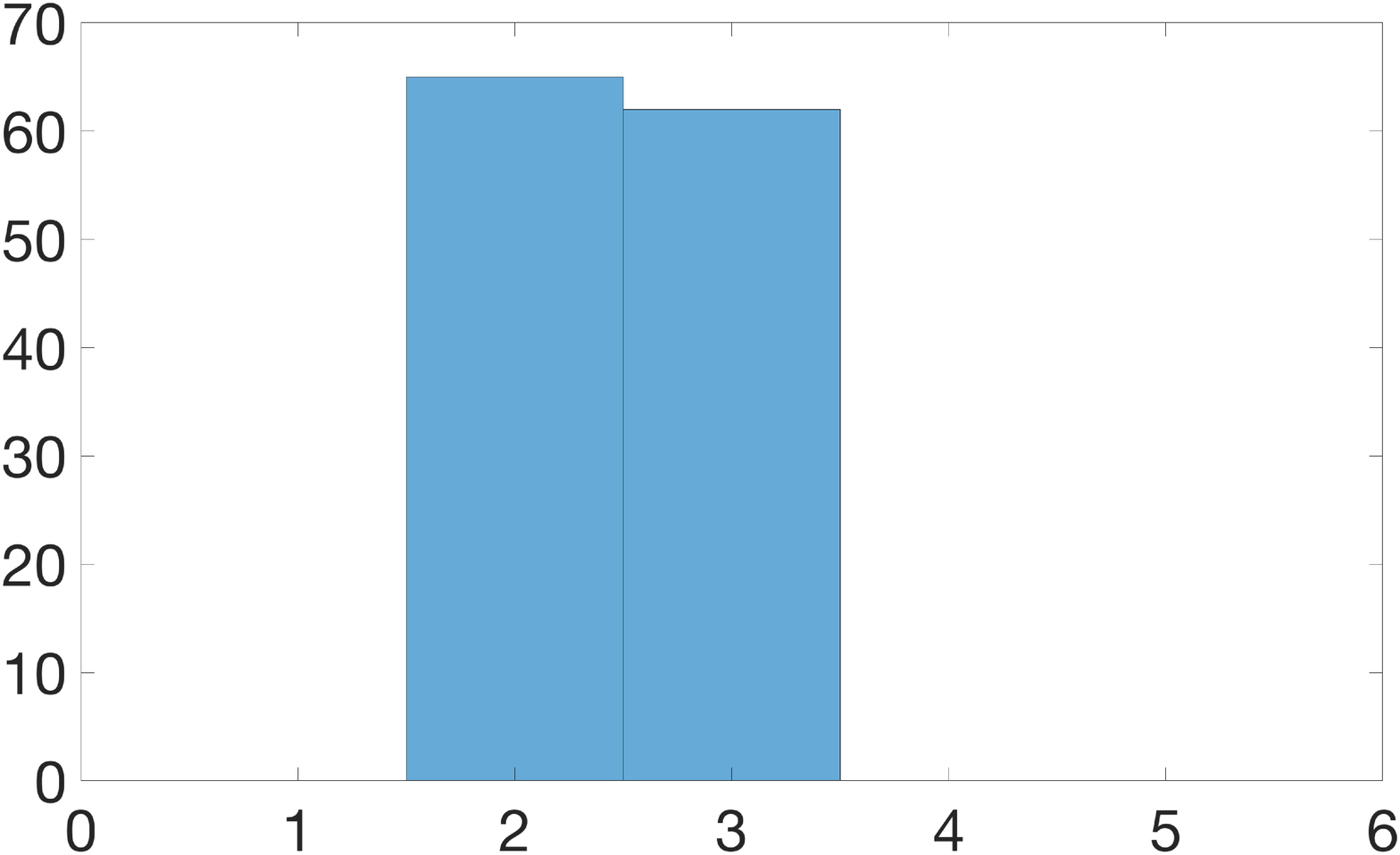}
\end{minipage} \\
$G_9$:
\begin{minipage}{65mm}
\includegraphics[width=60mm]{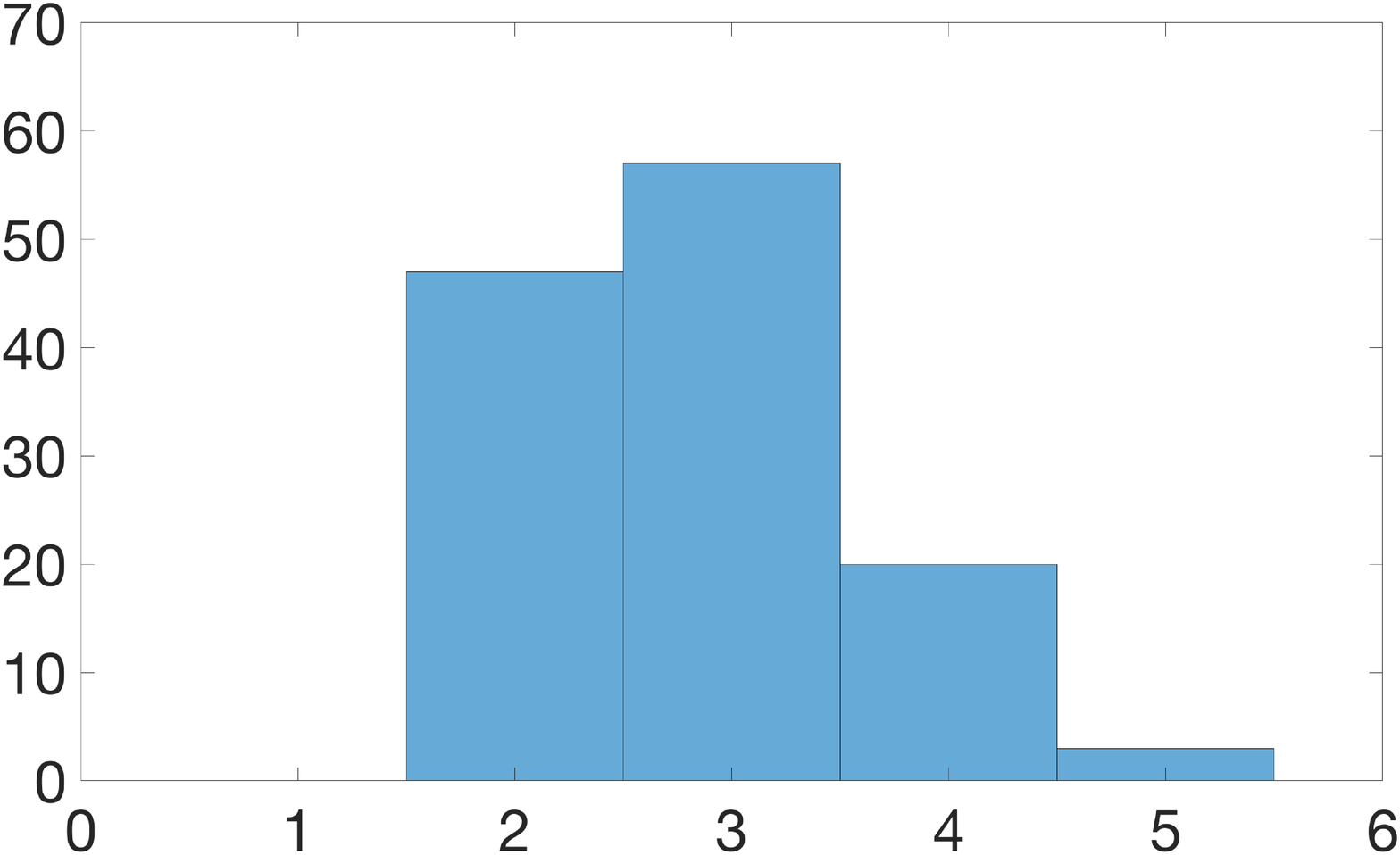}
\end{minipage}
\hspace{-1mm}$G_{10}$:
\begin{minipage}{65mm}
\includegraphics[width=60mm]{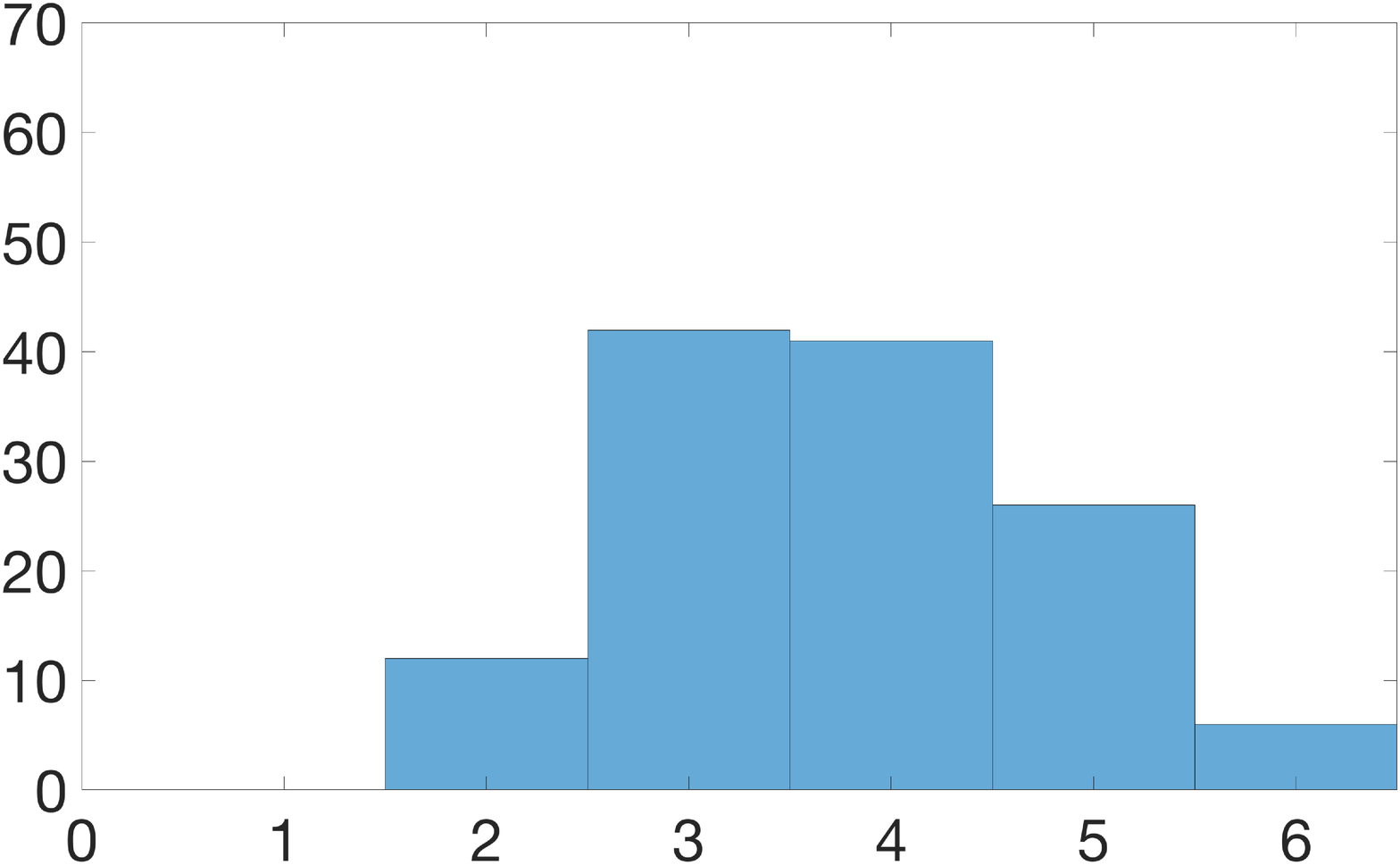}
\end{minipage} \\
\caption{\sf $N=127$: Degree distributions of the graphs $G_5$--$G_{10}$ obtained by Algorithm~\ref{algo1}.}
\label{fig:V127_A_deg_combined}
\end{figure}

\subsubsection{127-node networks}\label{127-nodes}

To understand the influence of the penalizing terms in the relaxed and discretized Problem~\eqref{model-3}, we implemented Algorithm~\ref{algo1} by solving first Problem~\eqref{model-3} using only the $\ell_1$-objective and ignoring all other terms which penalize a wide distribution of the degrees.  This produced a graph with 178 edges, $Q_L\approx 300$, and degrees ranging from 1 to 13.  This degree distribution is obviously not desirable.  

Including the terms in \eqref{model-3} penalizing wider degree distributions across the graph nodes, i.e., setting $c_2 = c_3 = 100$, Algorithm~\ref{algo1} has produced more desirable graphs, e.g., the graphs $G_6$ and $G_7$, the features $\E[\bar r]$, $\Var[\bar r]$ and the degree histograms of which are depicted in Figures~\ref{fig:V127_combined} and \ref{fig:V127_A_deg_combined}.  The generated graphs are still sparse and have much better spectral ratios $Q_L$.

As mentioned in Section \ref{tree}, we have also constructed a 127-node graph, namely $G_5$, using the tree--expander approach introduced and described in \cite{Taylor2020}.  All of $G_5$, $G_6$ and $G_7$ show somewhat similar/close characteristics in terms of the $\E[\bar r]$ and $\Var[\bar r]$ curves, although the graph featured in Figure~\ref{fig:V127_A_deg_combined} clearly has a much better degree distribution; its nodes have either degree 2 or degree 3.  A closer examination of the $\E[\bar r]$ and $\Var[\bar r]$ curves suggests that $G_7$ will synchronize more robustly than $G_5$ and $G_6$, since its $\E[\bar r]$ curve is above those of $G_5$ and $G_6$ and its $\Var[\bar r]$ curve is below those of $G_5$ and $G_6$, for $\sigma \ge 1.4$.

The tree--expander construction of a graph that is introduced and described in \cite{Taylor2020} allows/requires choices amongst many options; so it is possible to get quite different graphs by this construction.  A different graph constructed using an alternative sequence of choices is reported in \cite{Taylor2020}, the features of which can also be seen in Figure~\ref{fig:V127_combined}, labelled $G_8$.  This graph has better $\E[\bar r]$ and $\Var[\bar r]$ curves, in that the $\E[\bar r]$ curve approaches 1.0 relatively more quickly and the peak of $\Var[\bar r]$ is smaller as well as it dying down more quickly, compared with the previous three graphs discussed, namely $G_5$, $G_6$ and $G_7$.

A desirable/useful characteristic of a graph would presumably be the following: (i) $\E[\bar r]$ is almost ``S-shaped'' and approaches 1.0 rapidly and (ii) $\Var[\bar r]$ peaks early and dies down quickly.  Such a characteristic seems to be attainable by adding more edges, via an investigation using our optimization approach in Algorithm~\ref{algo1}.  In what follows, we consider solutions which have more edges, namely the graphs with 180 and 240 edges, respectively.

One might then argue that by increasing the number of edges the sparsity might be  compromised.  This would really depend on a particular application and what (number of edges) is practically meant by a sparse graph.  Define the {\em sparsity of an undirected graph} as the number of edges divided by $(N(N-1)/2)$.  In the previously discussed 127-node graphs the sparsity is 98.03\% ($G_5$ and $G_8$, 158 edges), 97.93\% ($G_6$, 166 edges), 97.91\% ($G_7$, 167 edges).  In the next two example graphs, the sparsity is 97.75\% (180 edges) and 97.00\% (240 edges), respectively.  In all of the 127-node examples presented in this paper, the sparsity might be regarded to be quite close to one another (depending on a particular application, of course), which in this case is around 97--98\%.

Figure~\ref{fig:V127_combined} features the optimized graph $G_9$ found by using Algorithm~\ref{algo1}, which has 180 edges, instead of the 158--167 edges the previous graphs have.  Although $G_9$ has more edges than the previous graphs, this may not necessarily be regarded to be compromising sparsity (as discussed above). Now, both of the $\E[\bar r]$ and $\Var[\bar r]$ curves have markedly better behaviour: The peak of $\Var[\bar r]$ is shifted to the left and $\E[\bar r]$ approaches 1.0 markedly more rapidly.  One might point to the rather wider range of degrees compared to the previous graphs; however, again depending on the application, this ``slightly'' wider range might be acceptable.

Figure~\ref{fig:V127_combined} also features the optimized graph $G_{10}$, again found by using Algorithm~\ref{algo1}, which has 240 edges.  Although 240 edges sounds to be much bigger than 160 edges, the sparsity of $G_{10}$ is worse only by 1 percentage point, which might again be acceptable in some practical situations. The behaviour of the graph shown by the $\E[\bar r]$ and $\Var[\bar r]$ curves is rather impressive:  The curve of $\E[\bar r]$ is almost S-shaped, settling very near 1.0 quite early.  The peak of $\Var[\bar r]$ has been pushed to the far left and $\Var[\bar r]$ dies down to almost zero around $\sigma = 1.0$ where the graph synchronize very well.  These point to a desirable robust synchrony of the graph.

\subsubsection{Computational effort}

\begin{table}[t]
\begin{center}
{\footnotesize
\begin{tabular}{lrrrr}
& \multicolumn{1}{c}{minimum\ \ \ \ \ } & \multicolumn{1}{c}{average\ \ \ \ \ } & \multicolumn{1}{c}{maximum\ \ \ \ \ } & \\
\multicolumn{1}{c}{$N$} & \multicolumn{1}{c}{CPU time [s]} & \multicolumn{1}{c}{CPU time [s]} & \multicolumn{1}{c}{CPU time [s]} & \multicolumn{1}{c}{Success rate} \\[1mm] \hline
\ \,20 & 3 & 5  & 8 & 17/30 (57\%) \\
\ \,30 & 10 & 18 & 54 & 23/30 (77\%) \\
\ \,50 & 64 & 139 & 365 & 14/30 (47\%) \\
127 & 2038 & 4419 & 7828 & 23/30 (77\%) \\[1mm]
\hline
 \end{tabular}}
\end{center}
\caption{\sf CPU times (in seconds) and success rates in solving model~\eqref{model-3} with various graph sizes~$N$.}
\label{CPU}
\end{table}

We have run AMPL--Ipopt suite to solve the optimization model~\eqref{model-3}, with randomly generated data, 30 times, with the graph sizes $N = 20$, 30, 50, and 127, on the computer the specifications of which we provide in Section~\ref{sec:implement}.  As can be seen from Table~\ref{CPU}, the computational suite was successful in about slightly more than half the time it was run.  As can also be seen, the elapsed computational time grew exponentially with the graph size, as expected.  Despite the exponential computational complexity, (locally) optimized graphs can be obtained in a reasonable amount of time.  For example, a locally optimal 127-node graph can be designed in about one to two hours' time.

Once $\widetilde{A}$ is found by solving model~\eqref{model-3} and an adjacency matrix $A$ is constructed by Algorithm~\ref{algo1}, the expected value $\E[\bar r]$ and the variance $\Var[\bar r]$ of $\bar r$ are computed and plotted against the coupling strength $\sigma$.  This computation also incurs a comparable CPU time, in addition to the CPU times reported in Table~\ref{CPU}.  As we mention in Section~\ref{sec:numexp}, we use the {\sc Matlab} solver {\tt ode15s}, effective for stiff ODEs, to solve ~\eqref{Kura}, which can take a long time depending on the randomly generated data.  In the case of $N=20$, the CPU times for the graphs $G_1$ and $G_4$ are about 6.2 and 4.3 minutes, respectively.  Given the fact that $G_4$ is more robust (although less sparse) than $G_1$, it is fair to say that these computations are expected to take shorter time for $G_1$ than those for $G_4$.  We observe a similar situation  for the $N=127$ case:  the CPU times for the graphs $G_5$ and $G_{10}$ are about 180 and 70 minutes, respectively.


\section{Conclusion and Future Work}
\label{conclusion}

We have proposed an optimization algorithm (Algorithm~\ref{algo1}) for designing sparse undirected networks, i.e., graphs, which synchronize well.  The algorithms involve solution of a discretized relaxed mathematical model which maximize sparsity and minimize the spread of degrees (in some sense) at the same time.  Synchrony of the network is posed as a constraint by means of the well-known Kuramoto system of coupled oscillators.  By means of carefully chosen examples, from small to large scale, we illustrate the working of our algorithm and optimization model.
The outcome is a method to generate
sparse well balanced graphs with good
synchronization from a random set
of frequencies and initial phases.
These graphs retain these properties when other frequency choices are used. In other words, in the resultant graphs there is not a strong correlation between the frequencies and degrees of the nodes (with
a cubic generated in one case)
underlying this insensitivity.

We have illustrated that our algorithm is successful in designing sparse synchronizing graphs with a relatively narrow distribution of degrees and small Laplacian eigenvalue ratio $Q_L$.  By allowing more edges in the graph, but without sacrificing the sparsity much, we managed to obtain impressive synchronization behaviour in terms of the expected order parameter curve $\E[\bar r]$ and its variance $\Var[\bar r]$.  These graphs synchronize in a robust manner in that for relatively small synchronizing values of $\bar r$ the variance of $\bar r$ under randomized data is small.

As pointed out in Remark~\ref{rem:vital}, Algorithm~\ref{algo1} might be used as a valuable tool in identifying the so-called vital nodes and the critical set of links of the optimized sparse networks.

We have not explicitly incorporated the spectral ratio $Q_L$ into our optimization model.  However, the numerical experiments have demonstrated that the additional terms involving the minimum and maximum degrees in the objective function promote smaller $Q_L$.  Having said this, a measure of $Q_L$ might still be directly/explicitly included in a future model by using the Rayleigh quotient and the representation properties of the Laplacian using quadratic functions, such as the one in \cite[Equation (14)]{Mohar}. Similar expressions that can be investigated can be found in \cite[Remark 4.2]{FazlyabDorflerPreciado}.

Another interesting feature of a future optimization model would be constraints imposed on the degrees which the nodes of a graph should have.  This would result in more targeted design of networks with certain desired load distributions.

\section*{Acknowledgments}
The authors would like to acknowledge valuable discussions with Subhra Dey at the initiation of this project and with Richard Taylor during its later stages.  This research was a collaboration under the auspices of the Modelling Complex Warfighting initiative between the Commonwealth of Australia (represented by the Defence Science and Technology Group) and the University of South Australia through a Defence Science Partnerships agreement.



\begin{thebibliography}{30}
\bibitem{Forger2017}
{D. B. Forger},
{\em Biological Clocks, Rhythms, and Oscillations},
(The MIT Press, Cambridge, Massachusetts, 2017)
\bibitem{Kuramoto84}
Y. Kuramoto,
\textsl{Chemical Oscillations, Waves, and Turbulence} (Springer, Berlin, 1984)

\bibitem{rodrigues_survey}
{F. A. Rodrigues, T. K. D. M. Peron, P. Ji, and J. Kurths}, 
{The Kuramoto model in complex networks}, Phys. Rep. {\bf 610}, 1--98 (2016) 


\bibitem{DekTay2013}
A. H. Dekker and R Taylor,
Synchronization properties of trees in the Kuramoto model,
SIAM J. Appl. Dyn. Syst. {\bf 12}, 596-–617 (2013)


\bibitem{RoggAey2004}
J. A. Rogge and D. Aeyels,
Stability of phase locking in a ring of unidirectionally coupled oscillators, 
J. Phys. A: Math. Gen. {\bf 37}, 11135--11148 (2004)
\bibitem{OchGor2010}
J Ochab, and P. F. Gora, 
Synchronization of coupled oscillators in a local one-dimensional Kuramoto model, 
Acta Physica Polonica B Proceedings Supplement, {\bf 3}, 453--462 (2010)
\bibitem{Ich2004}
T. Ichinomiya, 
Frequency synchronization in a random oscillator network, 
Phys. Rev. E {\bf 70}(2), 026116 (2004)
\bibitem{Hong2002}
H. Hong, M. Y. Choi, and Beom Jun Kim,
Synchronization on small-world networks,
Phys. Rev. E {\bf 65}, 026139 (2002)
\bibitem{Oh2007}
E. Oh, D.-S. Lee, B. Kahng, and D. Kim,
Synchronization transition of heterogeneously coupled oscillators on scale-free networks,
Phys. Rev. E {\bf 75}, 011104 (2007)
\bibitem{Dekk2007}
A. H. Dekker,
Studying organisational topology with simple computational models,
J. Artif. Soc. Simul. {\bf 10}, 6 (2007)
\bibitem{Kall2020}
A. C. Kalloniatis, T. A. McLennan-Smith, D. O. Roberts,
Modelling distributed decision-making in command and control using stochastic network synchronisation,
European Journal of Operational Research {\bf 284}, 588--603 (2020)
\bibitem{Brede2008}
M. Brede,
Local versus global synchronization in networks of non-identical Kuramoto oscillators, 
Eur. Phys. J. B. {\bf 62}, 87 (2008)
%
\bibitem{TanAoy2008}
T. Tanaka and T. Aoyagi,
Optimal weighted networks of phase oscillators for synchronization,
Phys. Rev. E. {\bf 78},  046210 (2008)
\bibitem{KellyGottwald2011} 
{D. Kelly and G. A. Gottwald}, 
{On the topology of synchrony optimized networks of a Kuramoto-model with non-identical oscillators}, 
Chaos {\bf 21}, 025110 (2011)
\bibitem{YanMik2012}
T. Yanagita and A. S. Mikhailov,
Design of oscillator networks with enhanced synchronization tolerance against noise,
Phys. Rev. E. {\bf 85}, 056206 (2012)
\bibitem{Fazylab2017}
M. Fazlyab, F. Doerfler and V. M. Preciado, 
Optimal network design for synchronization of coupled oscillators, 
Automatica, {\bf 84}, 181 (2017)
\bibitem{PecCar1998}
L. M. Pecora and T. L. Carroll,
Master stability functions for synchronizes coupled systems,
Phys. Rev. Lett. {\bf 80}, 2109 (1998)
\bibitem{BarPec2002}
M. Barahona and L. M. Pecora,
Synchronization in small world systems,
Phys. Rev. Lett. {\bf 89}, 054101 (2002)
\bibitem{Donetti2005}
L. Donetti, P. I. Hurtado, and M. A. Munoz,
Entangled networks, synchronization, and optimal network topology,
Phys. Rev. Lett. {\bf 95}, 188701 (2005)
\bibitem{Donetti2006}
L. Donetti, F. Neri, and M. A. Munoz,
Optimal network topologies: expanders, cages, Ramanujan
graphs, entangled networks and all that, 
J. Stat. Mech. P08007  (2006)
\bibitem{Estrada2010}
E. Estrada, S. Gago, and G. Caporossi,
Design of highly synchronizable and robust networks
Automatica,  {\bf 46}, 1835 (2010)
\bibitem{Taylor2020}
{R. Taylor, A. Kalloniatis, K. Hoek},
{Organisational hierarchy constructions with easy Kuramoto synchronisation}
J.Phys.A: Math and Theor, doi.10.1088/1751-8121/ab69a3
\bibitem{NocWri2006}
\newblock {J. Nocedal and S. Wright},
\newblock {\em Numerical Optimization}, (Springer, New York, 2006)
\bibitem{VosMau2006}
{G. Vossen and H. Maurer},
{On $L^1$-minimization in optimal control and applications to robotics}.
Optimal Control Applications and Methods. {\bf 27}, 301--321 (2006)
\bibitem{Kaya2010}
{C. Y. Kaya},
{Inexact restoration for Runge-Kutta discretization of optimal control problems},
SIAM J. Numer. Anal. {\bf 48}(4), 1492--1517 (2010)
\bibitem{KayMar2007}
{C. Y. Kaya and J. M. Mart{\'\i}nez},
{Euler discretization for inexact restoration and optimal control},
J. Optim. Theory Appl. {\bf 134}, 191--206 (2007)
\bibitem{AltKaySch2016}
{W. Alt, C. Y. Kaya, and C. Schneider},
{Dualization and discretization of linear-quadratic control problems with bang--bang solutions},
EURO J Comput. Optim. {\bf 4}, 47--77 (2016)
\bibitem{BanKay2013}
{N. Banihashemi and C. Y. Kaya},
{Inexact restoration for Euler discretization of box-constrained optimal control problems},
J. Optim. Theory Appl. {\bf 156}, 726--760 (2013)
\bibitem{BurKayMaj2014}
{R. S. Burachik, C. Y. Kaya, and S. N. Majeed},
{A duality approach for solving control-constrained linear-quadratic optimal control problems},
SIAM J. Control Optim. {\bf 52}, 1771--1782 (2014).


\bibitem{KayMau2014}
{C. Y. Kaya and H. Maurer}, 
{A numerical method for nonconvex multi-objective optimal control problems},
Comput. Optim. Appl. {\bf 57}(3), 685--702 (2014)
\bibitem{AMPL}
{R. Fourer, D. M. Gay, and B. W. Kernighan},
{\em AMPL: A Modeling Language for Mathematical Programming}, 2nd  ed. (Brooks/Cole Publishing Company/Cengage Learning, 2003)
\bibitem{WacBie2006}
{A. W\"achter and L. T. Biegler},
{On the implementation of a primal-dual interior point filter line
search algorithm for large-scale nonlinear programming}.
Math. Progr. {\bf 106}, 25--57 (2006)
\bibitem{VitalNodes}
{L. L\"u, D. Chen, X. L. Ren, Q. M. Zhang, Y. C. Zhang, T. Zhou},
{Vital nodes identification in complex networks},
Phys. Rep. {\bf 650}, 1--63 (2016)
\bibitem{adjacency_matrices}
{R. S. Burachik, A. C. Kalloniatis, and C. Y. Kaya}, Ancillary files for the preprint arXiv:2006.00428v1.
URL: https://arxiv.org/src/2006.00428v1/anc.  (2020)
\bibitem{Mohar}
{B. Mohar},
{\em Some applications of Laplace eigenvalues of graphs}, in Graph Symmetry: Algebraic
Methods and Applications, edited by G. Hahn and G. Sabidussi, NATO ASI Series C 497
(Kluwer, Dordrecht, Boston), 225--275 (1997)
\bibitem{FazlyabDorflerPreciado}
{M. Fazlyab, F. D\"orfler, and V. Preciado}, 
{Optimal network design for synchronization of coupled oscillators}, Automatica, {\bf 84}(C), 181--189 (2017)






\end{thebibliography}
\end{document}